\documentclass[aos,MSNbibl,seceqn,citesort,dvips]{arximspdf}
\usepackage{multirow}
\usepackage{dcolumn}

\doi{10.1214/11-AOS934}
\volume{39}
\issue{6}
\pubyear{2011}
\firstpage{3032}
\lastpage{3061}

\makeatletter

\newcommand{\csp}{\,}

\newcolumntype{d}[1]{D{.}{.}{#1}}

\newtheorem{theorem}{Theorem}[section]
\newtheorem{lem}{Lemma}[section]
\newproclaim{Remark}{Remark}

\makeatother

\begin{document}
\begin{frontmatter}

\title{A sieve M-theorem for bundled parameters in semiparametric
models, with application to the efficient estimation in a linear
model for censored~data}
\runtitle{M-theorem for bundled parameters}

\begin{aug}
\author[A]{\fnms{Ying} \snm{Ding}\thanksref{T1}\ead[label=e1]{yingding@umich.edu}}
\and
\author[A]{\fnms{Bin} \snm{Nan}\corref{}\thanksref{T1,T2}\ead[label=e2]{bnan@umich.edu}}
\runauthor{Y. Ding and B. Nan}
\affiliation{University of Michigan}
\address[A]{Department of Biostatistics \\
University of Michigan \\
1420 Washington Heights\\
Ann Arbor, Michigan 48109-2029\\
USA\\
\printead{e1}\\
\phantom{E-mail: }\printead*{e2}} 
\end{aug}

\thankstext{T1}{Supported in part by NSF Grant DMS-07-06700.}

\thankstext{T2}{Supported in part by NSF Grant DMS-10-07590
and NIH Grant R01-AG036802.}

\received{\smonth{6} \syear{2011}}
\revised{\smonth{9} \syear{2011}}

%
\begin{abstract}
In many semiparametric models that are parameterized by two types of
parameters---a Euclidean parameter of interest and an
infinite-dimensional nuisance parameter---the two parameters are
bundled together, that is, the nuisance parameter is an unknown
function that contains the parameter of interest as part of its
argument. For example, in a linear regression model for censored
survival data, the unspecified error distribution function involves the
regression coefficients. Motivated by developing an efficient
estimating method for the regression parameters, we propose a general
sieve M-theorem for bundled parameters and apply the theorem to
deriving the asymptotic theory for the sieve maximum likelihood
estimation in the linear regression model for censored survival data.
The numerical implementation of the proposed estimating method can be
achieved through the conventional gradient-based search algorithms such
as the Newton--Raphson algorithm. We show that the proposed estimator
is consistent and asymptotically normal and achieves the semiparametric
efficiency bound. Simulation studies demonstrate that the proposed
method performs well in practical settings and yields more efficient
estimates than existing estimating equation based methods. Illustration
with a real data example is also provided.
\end{abstract}

%

\end{frontmatter}

\section{Introduction} \label{secintro1}

In a semiparametric model that is parameterized by two types of
parameters---a finite-dimensional Euclidean parameter and an
infinite-dimensional parameter---oftentimes the infinite-dimensional
parameter is considered as a nuisance parameter, and the two parameters
are \textit{separated}. In many interesting statistical models, however,
the parameter of interest and the nuisance parameter are \textit{bundled
together}, a terminology used by Huang and Wellner
\cite{HuangWellner1997} when they reviewed the linear models under
interval censoring, which means that the infinite-dimensional parameter
is an unknown function of the parameter of interest. For example, in a
linear regression model for censored survival data, the unspecified
error distribution function, often treated as a nuisance parameter, is
a function of the regression coefficients. Other examples include the
single index model and the Cox regression model with an unspecified
link function.

There is a rich literature of asymptotic distributional theories for
M-estimation in a variety of semiparametric models with well-separated
parameters; see, for example, \cite
{HeEtal2010,Huang1996,HuangJian1999,Shen1997,WellnerZhang2007,ZhangHuangHua2010},
among
many others. Though many methodologies of M-estimation for bundled
parameters have been proposed in the literature, general asymptotic
distributional theories for such problems are still lacking. The only
estimation theories for bundled parameters we are aware of are the
sieve generalized method of moment of \cite{AiChen2003} and the
estimating equation approach of \cite{ChenEtal2003,NanEtal2009}.

In this article, we consider an extension of existing asymptotic
distributional theories
to accommodate situations where the estimation criteria are
parameterized with bundled parameters. The proposed theory has similar
flavor of Theorem 2 in~\cite{ChenEtal2003}, but they are different
because the latter requires an existing uniform consistent estimator of
the infinite-dimensional nuisance parameter with a convergence rate
faster than $n^{-1/4}$, which is then treated as a fixed function of
the parameter of interest in their estimating procedure, while we need
to simultaneously estimate both parameters through a sieve parameter
space; furthermore, their existing nuisance parameter estimator needs
to satisfy their condition (2.6), which is usually hard to verify when
its convergence rate is slower than $n^{-1/2}$. Our proposed theory is
general enough to cover a wide range of problems for bundled parameters
including the aforementioned single index model, the Cox model with
unknown link function and a linear model under different censoring
mechanisms. Rigorous proofs for each of the models, however, will take
lengthy derivations.
We only use the efficient estimation in the semiparametric linear
regression model with right censored data as an illustrative example
that motivates such a theoretical development and will present results
for other models elsewhere. Note that the considered example cannot be
directly put into the framework of restricted moments due to right
censoring, thus cannot be handled by the method of \cite{AiChen2003}.

Suppose that the failure time transformed by a known monotone
transformation is linearly related to a set of covariates, where the
failure time is subject to right censoring. Let $T_i$ denote the
transformed failure time and~$C_i$ denote the transformed censoring
time by the same transformation for subject $i$, $i=1,\ldots,n$. Let
$Y_i = \min(T_i,C_i)$ and $\Delta_i=I(T_i\leq C_i)$. Then the
semiparametric linear model we consider here can be written as
%
%
\begin{equation}\label{eqmodel1}
T_i = X_i'\beta_0 + e_{0,i},\qquad i=1,\ldots,n,
\end{equation}
where the errors $e_{0,i}$ are independent and identically distributed
(i.i.d.) with an unspecified distribution. When the failure time is
log-transformed, this model corresponds to the well-known accelerated
failure time model \cite{Kalbfleisch}. Here we assume that $(X_i,
C_i)$, $i=1,\ldots,n$, are i.i.d. and independent of $e_{0,i}$. This is
a common assumption for linear models with censored survival data,
which is particularly needed in \cite{RitovWellner1988} to derive the
efficient score function for~$\beta_0$. Such an assumption, however, is
stronger than necessary in the usual linear regression without
censoring, for which the error is only required to be uncorrelated with
covariates; see, for example, \cite{Chamberlain1987}. We also avoid
trivial transformations such as $\log(0)$ so that we always have $Y_i$'s
bounded from below.

The semiparametric linear regression model relates the failure time to
the covariates directly. It provides a straightforward interpretation
of the data and serves as an attractive alternative to the Cox model
\cite{Cox1972} in many applications.
Several estimators of the regression parameters have been proposed in
the literature since late 1970s, including the rank-based estimators
(see, e.g., \cite{Prentice,Wei1990,Tsiatis,Ying1993,Jin2003,Jin2006})
and the Buckley--James
estimator (see, e.g., \cite{BJ1979,Ritov,LaiYing1991}).
There are two major challenges in the estimation for such a~linear
model: (1) the estimating functions in the aforementioned methods are
discrete, leading to potential multiple solutions as well as numerical
difficulties; (2) none of the aforementioned methods is efficient.
Recently, Zeng and Lin \cite{ZengLin2007} developed a kernel-smoothed
profile likelihood estimating procedure for the accelerated failure
time model.
In this article, we consider a~sieve maximum likelihood approach for
model (\ref{eqmodel1}) for censored data. The proposed approach is
much intuitive, easy to implement numerically and asymptotically efficient.

It is easy to see that $T$ and $C$ are independent conditional on $X$
under the assumption $e_0 \perp(C,X)$. Hence the joint density
function of $Z=(Y,\Delta,X)$ can be written as
%
%
\begin{equation} \label{eqjdf}
f_{Y,\Delta,X}(y,\delta,x) = \lambda_0(y-x'\beta_0)^{\delta}\exp\{
-\Lambda_0(y-x'\beta_0)\}H(y,\delta,x),
\end{equation}
where $\Lambda_0(\cdot)$ is the true cumulative hazard function for the
error term $e_0$ and~$\lambda_0(\cdot)$ is its derivative. $H(y,\delta
,x)$ only depends on the conditional distribution of $C$ given $X$ and
the marginal distribution of $X$, and is free of $\beta_0$ and~$\lambda
_0$. To simplify the notation, we will ignore the factor $H$ from the
likelihood function. Then for i.i.d. observations $(Y_i,\Delta_i,X_i)$,
$i=1,\ldots,n$, from~(\ref{eqjdf}) we obtain the log
likelihood function for $\beta$ and $\lambda$ as
%
%
\begin{equation} \label{eqllk1}
l_n(\beta,\lambda)= n^{-1} \sum_{i=1}^{n} \biggl\{ \Delta_i \log\{
\lambda(Y_i-X_i'\beta)\}-\int I(Y_i \ge t) \lambda(t-X_i'\beta) \,dt
\biggr\}.\hspace*{-35pt}
\end{equation}
The log likelihood given in (\ref{eqllk1}) apparently is a
semiparametric model, where the argument of the nuisance parameter
$\lambda$ involves $\beta$; thus $\beta$ and $\lambda$ are bundled
parameters. To keep the positivity of $\lambda$, let $g(\cdot)=\log
\lambda(\cdot)$. Then the log likelihood function for $\beta$ and $g$,
using the counting process notation, can be written as
%
%
\begin{equation} \label{eqllk2}
l_n(\beta,g) = n^{-1}\sum_{i=1}^n \biggl\{\int g(t-X_i'\beta)
\,dN_i(t) - \int I(Y_i\geq t) e^{g(t-X_i'\beta)} \,dt\biggr\},\hspace*{-35pt}
\end{equation}
where $N_i(t)=\Delta_i I(Y_i \leq t)$ is the counting process for
subject $i$.

We propose a new approach by directly maximizing the log likelihood
function in a sieve space in which function $g(\cdot)$ is approximated
by B-splines. Numerically, the estimator can be easily obtained by the
Newton--Raphson algorithm or any gradient-based search algorithms. We
show that the proposed estimator is consistent and asymptotically
normal, and the limiting covariance matrix reaches the semiparametric
efficiency bound, which can be estimated either by inverting the
information matrix based on the efficient score function of the
regression parameters derived by \cite{RitovWellner1988}, or by
inverting the observed information matrix of all parameters, taking
into account that we are also estimating the nuisance parameters in the
sieve space for the log hazard function.

\section{The sieve M-theorem on the asymptotic normality of
semiparametric estimation for bundled parameters} \label{secgenThm}
In this section, we extend the general theorem introduced by
\cite{WellnerZhang2007}, which deals with the asymptotic normality of
semiparametric M-estimators of regression parameters when the
convergence rate of the estimator for nuisance parameters can be slower
than $n^{-1/2}$. In their theorem, the parameters of interest and the
nuisance parameters are assumed to be separated. We consider a more
general setting where the nuisance parameter can be a function of the
parameters of interest. The theorem is crucial in the proof of
asymptotic normality given in Theorem~\ref{thmnormality} for our
proposed estimators.

Some empirical process notation will be used from now on. We denote $Pf
= \int f(z) \,dP(z)$ and $\mathbb{P}_nf=n^{-1}\sum_{i=1}^n f(Z_i)$,
where $P$ is a probability measure, and $\mathbb{P}_n$ is an empirical
probability measure, and denote $\mathbb{G}_n f = n^{1/2}(\mathbb
{P}_n-P) f$. Given i.i.d. observations $Z_1,Z_2,\ldots,Z_n \in\mathcal
{Z}$, we estimate the unknown parameters $(\beta, \zeta(\cdot,\beta))$
by maximizing an objective function for $(\beta, \zeta(\cdot,\beta))$,
$n^{-1}\sum_{i=1}^n m(\beta,\zeta(\cdot,\beta);Z_i)=\mathbb{P}_n m(\beta
,\zeta(\cdot,\beta);Z)$, where $\beta$ is the parameter of interest,
and $\zeta(\cdot,\beta)$ is the nuisance parameter that can be a
function of $\beta$. Here ``$ \cdot$'' denotes the other arguments of
$\zeta$ besides $\beta$, which can be some components of $Z \in\mathcal
{Z}$. If the objective function $m$ is the log-likelihood function of a
single observation, then the estimator becomes the semiparametric
maximum likelihood estimator. Here we adopt similar notation in~\cite
{WellnerZhang2007}.

Let $\theta= (\beta,\zeta(\cdot,\beta))$, $\beta\in\mathcal{B}
\subset\mathbb{R}^d$ and $\zeta\in\mathcal{H}$, where $\mathcal{B}$ is
the parameter space of~$\beta$, and $\mathcal{H}$ is a class of
functions mapping from $\mathcal{Z} \times\mathcal{B}$ to $\mathbb
{R}$. Let $\Theta= \mathcal{B} \times\mathcal{H}$ be the parameter
space of $\theta$. Define a distance between $\theta_1, \theta_2 \in
\Theta$ by
\[
d(\theta_1, \theta_2) = \{|\beta_2 - \beta_1|^2 + \|\zeta_2(\cdot, \beta
_2) - \zeta_1(\cdot, \beta_1)\|^2\}^{1/2},
\]
where \mbox{$| \cdot|$} is the Euclidean distance, and $\| \cdot\|$ is some
norm. Let $\Theta_n$ be the sieve parameter space, a sequence of
increasing subsets of the parameter space $\Theta$ growing dense in
$\Theta$ as $n \to\infty$. We aim to find $\hat\theta_n \in\Theta_n$
such that $d(\hat\theta_n, \theta_0) = o_p(1)$ and $\hat\beta_n$ is
asymptotically normal.

For any fixed $\zeta(\cdot, \beta) \in\mathcal{H}$, let $\{\zeta_{\eta
}(\cdot, \beta)\dvtx\eta\mbox{ in a neighborhood of 0}\in\mathbb{R}\}$ be
a~smooth curve in $\mathcal{H}$ running through $\zeta(\cdot, \beta)$
at $\eta=0$, that is, $\zeta_{\eta}(\cdot, \beta)|_{\eta=0}=\zeta(\cdot,
\beta)$. Assume all $\zeta(\cdot,\beta) \in\mathcal{H}$ are at least
twice-differentiable with respect to $\beta$, and denote
\[
\mathbb{H}=\biggl\{h\dvtx h(\cdot, \beta)=\frac{\partial\zeta_{\eta}(\cdot,
\beta)}{\partial\eta}\bigg|_{\eta=0}, \zeta_{\eta} \in\mathcal{H},
\beta\in\mathcal{B}\biggr\}.
\]
Assume the objective function $m$ is twice Frechet differentiable.
Since for a small~$\delta$, we have $\zeta(\cdot,\beta+\delta)-\zeta
(\cdot,\beta) = \dot{\zeta}_{\beta}(\cdot,\beta)\delta+ o(\delta)$,
here $\dot\zeta_\beta(\cdot, \beta) = \partial\zeta(\cdot, \beta)
/\partial\beta$; then by the definition of functional derivatives it
follows that
\begin{eqnarray*}
&& \lim_{\delta\rightarrow0} \frac{1}{\delta} \bigl\{m\bigl(\beta,\zeta
(\cdot,\beta+\delta);z\bigr)-m(\beta,\zeta(\cdot,\beta);z) \bigr\} \\
&&\qquad = \lim_{\delta\rightarrow0} \frac{1}{\delta} \bigl\{m\bigl(\beta
,\zeta(\cdot,\beta)+\dot{\zeta}_{\beta}(\cdot,\beta)\delta+o(\delta);z\bigr)
\\
&&\hspace*{50pt}\qquad\quad{} - m\bigl(\beta,\zeta(\cdot,\beta)+\dot{\zeta}_{\beta
}(\cdot,\beta)\delta;z\bigr)\bigr\} \\
&&\qquad\quad{} + \lim_{\delta\rightarrow0} \frac{1}{\delta}
\bigl\{m\bigl(\beta,\zeta(\cdot,\beta)+\dot{\zeta}_{\beta}(\cdot,\beta)\delta;z\bigr)
- m(\beta,\zeta(\cdot,\beta);z)\bigr\} \\
&&\qquad = \lim_{\delta\rightarrow0} \dot{m}_2\bigl(\beta,\zeta(\cdot
,\beta)+\dot{\zeta}_{\beta}(\cdot,\beta)\delta;z\bigr)[o(\delta)/\delta] \\
&&\qquad\quad{} + \dot{m}_2(\beta,\zeta(\cdot,\beta);z)[\dot{\zeta
}_{\beta}(\cdot,\beta)] \\
&&\qquad = \dot{m}_2(\beta,\zeta(\cdot,\beta);z)[\dot{\zeta}_{\beta
}(\cdot,\beta)],
\end{eqnarray*}
where the subscript 2 indicates that the derivatives are taken with
respect to the second argument of the function. The last equality holds because
\[
\lim_{\delta\rightarrow0} \dot{m}_2\bigl(\beta,\zeta(\cdot,\beta)+\dot
{\zeta}_{\beta}(\cdot,\beta)\delta;z\bigr)[o(\delta)/\delta]=0.
\]
Similarly we have
\begin{eqnarray*}
&& \lim_{\delta\rightarrow0}\frac{1}{\delta}\bigl\{\dot{m}_2\bigl(\beta
,\zeta(\cdot,\beta+\delta);z\bigr)[h(\cdot,\beta)] -\dot{m}_2(\beta,\zeta
(\cdot,\beta);z)[h(\cdot,\beta)] \bigr\} \\
&&\qquad = \ddot{m}_{22}(\beta,\zeta(\cdot,\beta);z)[h(\cdot,\beta
),\dot{\zeta}_{\beta}(\cdot,\beta)]
\end{eqnarray*}
and
\begin{eqnarray*}
&& \lim_{\delta\rightarrow0} \frac{1}{\delta} \{\dot{m}_2(\beta
,\zeta(\cdot,\beta);z)[h(\cdot,\beta+\delta)]
- \dot{m}_2(\beta,\zeta(\cdot,\beta);z)[h(\cdot,\beta)] \} \\
&&\qquad = \dot{m}_2(\beta,\zeta(\cdot,\beta);z)[\dot{h}_{\beta}(\cdot
,\beta)].
\end{eqnarray*}
Thus according to the chain rule of the functional derivatives, we
have
\begin{eqnarray*}
\dot{m}_{\beta}(\beta,\zeta(\cdot,\beta);z) &=& \frac{\partial m(\beta
,\zeta(\cdot,\beta);z)}{\partial\beta} \\
&=& \dot{m}_1(\beta,\zeta(\cdot,\beta);z) + \dot{m}_2(\beta,\zeta(\cdot
,\beta);z)[\dot{\zeta}_{\beta}(\cdot,\beta)],\hspace*{-22pt}\\
\dot{m}_{\zeta}(\beta,\zeta(\cdot,\beta);z)[h] &=& \frac{\partial
m(\beta,(\zeta+\eta h)(\cdot,\beta);z)}{\partial\eta}\bigg|_{\eta=0}
\\
&=& \dot{m}_2(\beta,\zeta(\cdot,\beta);z)[h(\cdot,\beta)],
\end{eqnarray*}
\begin{eqnarray*}
\ddot{m}_{\beta\beta}(\beta,\zeta(\cdot,\beta);z) &=& \frac{\partial^2
m(\beta,\zeta(\cdot,\beta);z)}{\partial\beta\,\partial\beta'} = \frac
{\partial\dot{m}_\beta(\beta,\zeta(\cdot,\beta);z)}{\partial\beta'} \\
&=& \ddot{m}_{11}(\beta,\zeta(\cdot,\beta);z)+ \ddot{m}_{12}(\beta;\zeta
(\cdot,\beta);z)[\dot{\zeta}_{\beta}(\cdot,\beta)]\\
&&{} + \ddot{m}_{21}(\beta,\zeta(\cdot,\beta);z)[\dot{\zeta}_{\beta
}(\cdot,\beta)] \\
&&{} + \ddot{m}_{22}(\beta,\zeta(\cdot,\beta);z)[\dot{\zeta}_{\beta
}(\cdot,\beta),\dot{\zeta}_{\beta}(\cdot,\beta)] \\
&&{} + \dot{m}_{2}(\beta,\zeta(\cdot,\beta);z)[\ddot{\zeta}_{\beta
\beta}(\cdot,\beta)], \\
\ddot{m}_{\beta\zeta}(\beta,\zeta(\cdot,\beta);z)[h] &=& \frac
{\partial\dot{m}_{\beta}(\beta,(\zeta+\eta h)(\cdot,\beta
);z)}{\partial\eta} \bigg\vert_{\eta=0} \\
&=& \ddot{m}_{12}(\beta,\zeta(\cdot,\beta);z)[h(\cdot,\beta)] \\
&&{} + \ddot{m}_{22}(\beta,\zeta(\cdot,\beta);z)[\dot{\zeta}_{\beta
}(\cdot,\beta),h(\cdot,\beta)] \\
&&{} + \dot{m}_{2}(\beta,\zeta(\cdot,\beta);z)[\dot{h}_\beta(\cdot
,\beta)], \\
\ddot{m}_{\zeta\beta}(\beta,\zeta(\cdot,\beta);z)[h] &=& \frac
{\partial\dot{m}_2(\beta,\zeta(\cdot,\beta);z)[h(\cdot,\beta
)]}{\partial\beta} \\
&=& \ddot{m}_{21}(\beta,\zeta(\cdot,\beta);z)[h(\cdot,\beta)] \\
&&{} + \ddot{m}_{22}(\beta,\zeta(\cdot,\beta);z)[h(\cdot,\beta
),\dot{\zeta}_{\beta}(\cdot,\beta)] \\
&&{} + \dot{m}_{2}(\beta,\zeta(\cdot,\beta);z)[\dot{h}_\beta(\cdot
,\beta)], \\
\ddot{m}_{\zeta\zeta}(\beta,\zeta(\cdot,\beta);z)[h_1,h_2] &=&
\ddot{m}_{22}(\beta,\zeta(\cdot,\beta);z)[h_1(\cdot,\beta),h_2(\cdot
,\beta)].
\end{eqnarray*}

As noted before, the subscript 1 or 2 in the derivatives indicates that
the derivatives are taken with respect to the first or the second
argument of the function, and $h$ inside the square brackets is a
function denoting the direction of the functional derivative with
respect to $\zeta$. Note that for the second derivatives $\ddot
{m}_{\beta\zeta}$ and $\ddot{m}_{\zeta\beta}$, we implicitly require
the direction $h$ to be a differentiable function with respect to $\beta
$. It is easily seen that when $\zeta$ is free of $\beta$, all the
above derivatives reduce to that in~\cite{WellnerZhang2007}. Following
\cite{WellnerZhang2007}, we also define
\begin{eqnarray*}
\dot{S}_{\beta}(\beta,\zeta(\cdot,\beta)) &=& P\dot{m}_{\beta}(\beta
,\zeta(\cdot,\beta);Z), \\ \dot{S}_{\zeta}(\beta,\zeta(\cdot,\beta))[h]
&=& P\dot{m}_{\zeta}(\beta,\zeta(\cdot,\beta);Z)[h], \\
\dot{S}_{\beta, n}(\beta,\zeta(\cdot,\beta)) &=& \mathbb{P}_n \dot
{m}_{\beta}(\beta,\zeta(\cdot,\beta);Z), \\
\dot{S}_{\zeta, n}(\beta,\zeta(\cdot,\beta))[h] &=& \mathbb{P}_n \dot
{m}_{\zeta}(\beta,\zeta(\cdot,\beta);Z)[h], \\
\ddot{S}_{\beta\beta}(\beta,\zeta(\cdot,\beta)) &=& P\ddot{m}_{\beta
\beta}(\beta,\zeta(\cdot,\beta);Z), \\
\ddot{S}_{\zeta\zeta}(\beta,\zeta(\cdot,\beta))[h,h] &=& P\ddot
{m}_{\zeta\zeta}(\beta,\zeta(\cdot,\beta);Z)[h,h]
\end{eqnarray*}
and
\[
\ddot{S}_{\beta\zeta}(\beta,\zeta(\cdot,\beta))[h] = \ddot{S}_{\zeta
\beta}'(\beta,\zeta(\cdot,\beta))[h]
= P \ddot{m}_{\beta\zeta}(\beta, \zeta(\cdot,\beta);Z)[h].
\]
Furthermore, for $\mathbf{h}=(h_1,h_2,\ldots,h_d)' \in\mathbb{H}^d$,
we denote
\begin{eqnarray*}
\dot{m}_{\zeta}(\beta,\zeta(\cdot,\beta);z)[\mathbf{h}]
&=& (\dot{m}_{\zeta}(\beta,\zeta(\cdot,\beta);z)[h_1],\ldots, \dot
{m}_{\zeta}(\beta,\zeta(\cdot,\beta);z)[h_d])', \\
\ddot{m}_{\beta\zeta}(\beta,\zeta(\cdot,\beta
);z)[\mathbf{h}] &=& (\ddot{m}_{\beta\zeta}(\beta,\zeta(\cdot,\beta
);z)[h_1],\ldots, \ddot{m}_{\beta\zeta}(\beta,\zeta(\cdot,\beta
);z)[h_d]), \\
\ddot{m}_{\zeta\beta}(\beta,\zeta(\cdot,\beta
);z)[\mathbf{h}] &=& (\ddot{m}_{\zeta\beta}(\beta,\zeta(\cdot,\beta
);z)[h_1],\ldots, \ddot{m}_{\zeta\beta}(\beta,\zeta(\cdot,\beta
);z)[h_d])', \\
\ddot{m}_{\zeta\zeta}(\beta,\zeta(\cdot,\beta
);z)[\mathbf{h},h] &=& (\ddot{m}_{\zeta\zeta}(\beta,\zeta(\cdot,\beta
);z)[h_1,h],\ldots,\\
&&\hspace*{26.2pt}\ddot{m}_{\zeta\zeta}(\beta,\zeta(\cdot,\beta);z)[h_d,h])'
\end{eqnarray*}
and define correspondingly
\begin{eqnarray*}
\dot{S}_{\zeta}(\beta,\zeta(\cdot,\beta))[\mathbf{h}] &=& P \dot
{m}_{\zeta}(\beta,\zeta(\cdot,\beta);Z)[\mathbf{h}], \\
\dot{S}_{\zeta, n}(\beta,\zeta(\cdot,\beta))[\mathbf{h}] &=& \mathbb
{P}_n \dot{m}_{\zeta}(\beta,\zeta(\cdot,\beta);Z)[\mathbf{h}], \\
\ddot{S}_{\beta\zeta}(\beta,\zeta(\cdot,\beta))[\mathbf{h}] &=& P \ddot
{m}_{\beta\zeta}(\beta,\zeta(\cdot,\beta);Z)[\mathbf{h}], \\
\ddot{S}_{\zeta\beta} (\beta,\zeta(\cdot,\beta))[\mathbf{h}] &=& P
\ddot{m}_{\zeta\beta}(\beta,\zeta(\cdot,\beta);Z)[\mathbf{h}], \\
\ddot{S}_{\zeta\zeta}(\beta,\zeta(\cdot,\beta))[\mathbf{h},h] &=& P
\ddot{m}_{\zeta\zeta}(\beta,\zeta(\cdot,\beta);Z)[\mathbf{h},h].
\end{eqnarray*}

To obtain the asymptotic normality result for the sieve M-estimator
$\hat{\beta}_n$, the assumptions we will make in the following look
similar to those in \cite{WellnerZhang2007}, but all the derivatives
with respect to $\beta$ involve the chain rule and hence are more
complicated, which is the key difference to \cite{WellnerZhang2007}.
Additionally, we focus on sieve estimators in the sieve parameter
space. We list the following assumptions:
\begin{longlist}[(A2)]
\item[(A1)] (Rate of convergence) For an estimator $\hat\theta_n =
(\hat\beta_n, \hat\zeta_n(\cdot,\hat\beta_n)) \in\Theta_n$ and the
true parameter $\theta_0=(\beta_0, \zeta_0(\cdot,\beta_0))\in\Theta$,
$d(\hat\theta_n, \theta_0) = O_p(n^{-\xi})$ for some \mbox{$\xi>0$}.
\item[(A2)] $\dot{S}_\beta(\beta_0,\zeta_0(\cdot,\beta_0))=0$ and $\dot
{S}_\zeta(\beta_0,\zeta_0(\cdot,\beta_0))[h]=0$ for all $h \in\mathbb{H}$.
\item[(A3)] (Positive information) There exists an $\mathbf
{h}^*=(h_1^*,\ldots,h_d^*)'$, where\break \mbox{$h_j^* \in\mathbb{H}$} for
$j=1,\ldots,d$, such that
\[
\ddot{S}_{\beta\zeta}(\beta_0,\zeta_0(\cdot,\beta_0))[h]-\ddot
{S}_{\zeta\zeta}(\beta_0,\zeta_0(\cdot,\beta_0))[\mathbf{h}^*,h] = 0
\]
for all $h \in\mathbb{H}$. Furthermore, the matrix
\begin{eqnarray*}
A &=& -\ddot{S}_{\beta\beta}(\beta_0,\zeta_0(\cdot,\beta_0)) + \ddot
{S}_{\zeta\beta}(\beta_0,\zeta_0(\cdot,\beta_0))[\mathbf{h}^*] \\
&=& - P\{\ddot{m}_{\beta\beta}(\beta_0,\zeta_0(\cdot,\beta
_0);Z)-\ddot{m}_{\zeta\beta}(\beta_0,\zeta_0(\cdot,\beta_0);Z)[\mathbf
{h}^*]\}
\end{eqnarray*}
is nonsingular.
\item[(A4)] The estimator $(\hat{\beta}_n,\hat{\zeta}_n(\cdot, \hat
\beta_n))$ satisfies
\begin{eqnarray*}
\dot{S}_{\beta, n}(\hat{\beta}_n,\hat{\zeta}_n(\cdot,\hat{\beta}_n))=
o_p(n^{-1/2}) \quad\mbox{and}\quad \dot{S}_{\zeta,n}(\hat{\beta}_n,\hat{\zeta
}_n(\cdot,\hat{\beta}_n))[\mathbf{h}^*]
= o_p(n^{-1/2}).
\end{eqnarray*}

\item[(A5)] (Stochastic equicontinuity) For some $C>0$,
\begin{eqnarray*}
&& \sup_{d(\theta,\theta_0) \leq Cn^{-\xi}, \theta\in\Theta_n}
\bigl\vert\sqrt{n}(\dot{S}_{\beta, n}-\dot{S}_\beta)(\beta,\zeta(\cdot,\beta
)) \\
&&\hspace*{55pt}\qquad{} - \sqrt{n}(\dot{S}_{\beta, n}-\dot
{S}_\beta)(\beta_0,\zeta_0(\cdot,\beta_0))\bigr\vert= o_p(1)
\end{eqnarray*}
and
\begin{eqnarray*}
&& \sup_{d(\theta,\theta_0) \leq Cn^{-\xi}, \theta
\in\Theta_n} \bigl\vert\sqrt{n}(\dot{S}_{\zeta, n}-\dot{S}_\zeta)(\beta
,\zeta(\cdot,\beta))[\mathbf{h}^*(\cdot,\beta)] \\
&&\hspace*{54pt}\qquad{} - \sqrt{n}(\dot
{S}_{\zeta, n}-\dot{S}_\zeta)(\beta_0,\zeta_0(\cdot,\beta_0))[\mathbf
{h}^*(\cdot,\beta_0)]\bigr\vert=o_p(1).
\end{eqnarray*}
\item[(A6)] (Smoothness of the model) For some $\alpha>1$ satisfying
$\alpha\xi>1/2$, and for $\theta$ in a neighborhood of $\theta_0\dvtx\{
\theta\dvtx d(\theta, \theta_0) \leq Cn^{-\xi}, \theta\in\Theta_n \}$,
\begin{eqnarray*}
&& \vert\dot{S}_\beta(\beta,\zeta(\cdot,\beta
))-\dot{S}_\beta(\beta_0,\zeta_0(\cdot,\beta_0))
- \ddot{S}_{\beta\beta}(\beta_0,\zeta_0(\cdot,\beta_0))(\beta-\beta
_0) \\
&&\hspace*{81pt}\qquad{} - \ddot{S}_{\beta\zeta}(\beta_0,\zeta_0(\cdot,\beta
_0))[\zeta(\cdot,\beta)-\zeta_0(\cdot,\beta_0)]\vert\\
&&\qquad= O(d^\alpha
(\theta,\theta_0))
\end{eqnarray*}
and
\begin{eqnarray*}
&& \vert\dot{S}_\zeta(\beta,\zeta(\cdot,\beta
))[\mathbf{h}^*(\cdot,\beta)] -\dot{S}_\zeta(\beta_0,\zeta_0(\cdot,\beta
_0))[\mathbf{h}^*(\cdot,\beta_0)] \\
&&\quad{} - \ddot{S}_{\zeta\beta}(\beta_0,\zeta_0(\cdot,\beta_0))[\mathbf
{h}^*(\cdot,\beta_0)](\beta-\beta_0) \\
&&\hspace*{5pt}\quad{} - \ddot{S}_{\zeta\zeta
}(\beta_0,\zeta_0(\cdot,\beta_0))[\mathbf{h}^*(\cdot,\beta_0),\zeta
(\cdot,\beta)-\zeta_0(\cdot,\beta_0)] \vert\\
&&\hspace*{5pt}\qquad= O(d^\alpha(\theta
,\theta_0)).
\end{eqnarray*}
\end{longlist}
Note that $\xi$ in (A1) depends on the entropy of the sieve parameter
space for~$\zeta$ and cannot be arbitrarily small; it is controlled by
the smoothness of the model in (A6). The convergence rate in (A1) needs
to be achieved prior to obtaining asymptotic normality. Assumption (A2) is a
common assumption for the maximum likelihood estimation and usually
holds. The direction $\mathbf{h}^*$ in (A3) may be found through the
equation in (A3). It is the least favorable direction when $m$ is the
likelihood function. Assumptions (A4) and (A5) are usually verified either by the
Donsker property or the maximal inequality of~\cite{vanWellner1996}.
Assumption (A6) can be obtained by a Taylor expansion. The following theorem is an
extension to Theorem 6.1 in \cite{WellnerZhang2007} when the
infinite-dimensional parameter $\zeta$ is a function of the
finite-dimensional parameter~$\beta$.
%
%
\begin{theorem}\label{thmgenThm}
Suppose that assumptions \textup{(A1)--(A6)} hold. Then
\begin{eqnarray*}
&&\sqrt{n}(\hat{\beta}_n-\beta_0) = A^{-1}\sqrt{n} \mathbb{P}_n
m^*(\beta_0,\zeta_0(\cdot,\beta_0);Z)+o_p(1) \\
&&\hspace*{61.64pt}\rightarrow_d N(0,A^{-1}B(A^{-1})'),
\end{eqnarray*}
where
\begin{eqnarray*}
m^*(\beta_0,\zeta_0(\cdot,\beta_0);z) &=& \dot{m}_\beta(\beta_0,\zeta
_0(\cdot,\beta_0);z)- \dot{m}_\zeta(\beta_0,\zeta_0(\cdot,\beta
_0);z)[\mathbf{h}^*], \\
B &=& P \{ m^*(\beta_0,\zeta_0(\cdot,\beta_0);Z)^{\otimes2}\},
\end{eqnarray*}
and $A$ is given in assumption \textup{(A3)}. Here $a^{\otimes2} = aa'$.
\end{theorem}
\begin{pf}
The proof follows similarly along the proof of Theorem 6.1 in~\cite
{WellnerZhang2007}. Assumptions (A1) and (A5) yield
\[
\sqrt{n}(\dot{S}_{\beta, n}-\dot{S}_\beta)(\hat{\beta}_n,\hat{\zeta
}_n(\cdot,\hat{\beta}_n))-\sqrt{n}(\dot{S}_{\beta, n}-\dot{S}_\beta
)(\beta_0,\zeta_0(\cdot,\beta_0))=o_p(1).
\]
Since $\dot{S}_{\beta, n}(\hat{\beta}_n,\hat{\zeta}_n(\cdot,\hat{\beta
}_n))=o_p(n^{-1/2})$ by (A4) and $\dot{S}_{\beta}(\beta_0,\zeta_0(\cdot
,\beta_0))=0$ by~(A2), we have
\[
\sqrt{n}\dot{S}_\beta(\hat{\beta}_n,\hat{\zeta}_n(\cdot,\hat{\beta
}_n))+ \sqrt{n}\dot{S}_{\beta, n}(\beta_0,\zeta_0(\cdot,\beta_0)) =o_p(1).
\]
Similarly,
\[
\sqrt{n}\dot{S}_\zeta(\hat{\beta}_n,\hat{\zeta}_n(\cdot,\hat{\beta
}_n))[\mathbf{h}^*(\cdot,\hat{\beta}_n)] + \sqrt{n}\dot{S}_{\zeta,
n}(\beta_0,\zeta_0(\cdot,\beta_0))[\mathbf{h}^*(\cdot,\beta_0)] =o_p(1).
\]
Combining these equalities and assumption (A6) yields
%
%
\begin{eqnarray}\label{eqgenThm1}
&& \ddot{S}_{\beta\beta}(\beta_0,\zeta_0(\cdot,\beta
_0))(\hat{\beta}_n-\beta_0)+\ddot{S}_{\beta\zeta}(\beta_0,\zeta_0(\cdot
,\beta_0))[\hat{\zeta}_n(\cdot,\hat{\beta}_n)-\zeta_0(\cdot,\beta_0)]
\nonumber\\
&&\quad{} + \dot{S}_{\beta, n}(\beta_0,\zeta_0(\cdot,\beta_0)) + O(d^\alpha
(\hat\theta_n, \theta_0)) \\
&&\qquad = o_p(n^{-1/2})\nonumber
\end{eqnarray}
and
%
%
\begin{eqnarray} \label{eqgenThm2}
&& \ddot{S}_{\zeta\beta}(\beta_0,\zeta_0(\cdot,\beta
_0))[\mathbf{h}^*(\cdot,\beta_0)](\hat{\beta}_n-\beta_0)\nonumber\\
&&\quad{}+\ddot
{S}_{\zeta\zeta}(\beta_0,\zeta_0(\cdot,\beta_0))[\mathbf{h}^*(\cdot
,\beta_0),  \hat{\zeta}_n(\cdot,\hat{\beta}_n) -
\zeta_0(\cdot,\beta_0)]\nonumber\\[-8pt]\\[-8pt]
&&\quad{} + \dot
{S}_{\zeta, n}(\beta_0,\zeta_0(\cdot,\beta_0))[\mathbf{h}^*(\cdot,\beta
_0)] + O(d^\alpha(\hat\theta_n, \theta_0)) \nonumber\\
&&\qquad = o_p(n^{-1/2}).\nonumber
\end{eqnarray}
Since $\alpha>1$ with $\alpha\xi>1/2$, the rate of convergence
assumption (A1) implies $\sqrt{n} O(d^\alpha(\hat\theta_n, \theta_0)) =
o_p(1)$, then (\ref{eqgenThm1}) and (\ref{eqgenThm2}) together with (A3) yields
\begin{eqnarray*}
&& \bigl(\ddot{S}_{\beta\beta}(\beta_0,\zeta_0(\cdot,\beta_0))-\ddot
{S}_{\zeta\beta}(\beta_0,\zeta_0(\cdot,\beta_0))[\mathbf{h}^*(\cdot
,\beta_0)]\bigr)(\hat{\beta}_n-\beta_0) \\
&&\qquad = - \bigl(\dot{S}_{\beta, n}(\beta_0,\zeta_0(\cdot,\beta_0))-\dot
{S}_{\zeta, n}(\beta_0,\zeta_0(\cdot,\beta_0))[\mathbf{h}^*(\cdot,\beta
_0)]\bigr)+o_p(n^{-1/2}),
\end{eqnarray*}
that is,
\[
- A(\hat{\beta}_n-\beta_0) = -\mathbb{P}_n m^*(\beta_0,\zeta_0(\cdot
,\beta_0);Z)+o_p(n^{-1/2}).
\]
This yields
\begin{eqnarray*}
&&\sqrt{n}(\hat{\beta}_n-\beta_0) = A^{-1}\sqrt{n}\mathbb{P}_n m^*(\beta
_0,\zeta_0(\cdot,\beta_0);Z) + o_p(1) \\
&&\hspace*{60.64pt}\rightarrow_d N(0,A^{-1}B(A^{-1})').
\end{eqnarray*}
\upqed\end{pf}

\section{Back to the linear model: The sieve maximum likelihood
estimation} \label{secspline}
By taking logarithm to the positive function $\lambda(\cdot)$ in (\ref
{eqllk1}), the function $g(\cdot)$ in (\ref{eqllk2}) is no longer
restricted to be positive, which eases the estimation. We now describe
the spline-based sieve maximum likelihood estimation for model (\ref
{eqmodel1}). Under the regularity conditions (C.1)--(C.3) stated in
Section~\ref{secasymp}, we know that the observed residual times $\{
Y_i-X_i'\beta\dvtx\beta\in\mathcal{B}, i=1,\ldots,n\}$ are confined in
some finite interval. Let $[a,b]$ be an interval of interest, where
$-\infty<a<b<\infty$. Let $T_{K_n}=\{t_1,\ldots,t_{K_n}\}$ be a set of
partition points of $[a,b]$ with $K_n = O(n^{\nu})$ and ${\max_{1\leq j
\leq K_n+1}}|t_j-t_{j-1}|=O(n^{-\nu})$ for some constant $\nu\in(0, 1/2)$.
Let $\mathcal{S}_n(T_{K_n},K_n,p)$ be the space of polynomial splines
of order $p \geq1$ defined in \cite{Schumaker}, Definition 4.1.
According to Schumaker (\cite{Schumaker}, Corollary~4.10), there exists
a set of B-spline basis functions $\{B_j, 1 \leq j \leq q_n\}$ with
$q_n=K_n+p$ such that for any $s \in\mathcal{S}_n(T_{K_n},K_n,p)$, we
can write
%
%
\begin{equation} \label{eqgfun}
s(t) = \sum_{j=1}^{q_n} \gamma_j B_j(t),
\end{equation}
where we follow \cite{ShenWong1994} by requiring ${\max_{j=1, \ldots,
q_n}}|\gamma_j| \leq c_n$ that is allowed to grow with $n$ slowly enough.

Let $\gamma=(\gamma_1,\ldots, \gamma_ {q_n})'$. Under suitable
smoothness assumptions, $g_0(\cdot)=\log\lambda_0(\cdot)$ can be well
approximated by some function in $\mathcal{S}_n(T_{K_n},K_n,p)$.
Therefore,\vadjust{\goodbreak} we seek a member of $\mathcal{S}_n(T_{K_n},K_n,p)$ together
with a value of $\beta\in\mathcal{B}$ that maximizes the log
likelihood function. Specifically, let $\hat{\theta}_n = (\hat{\beta
}_n,\hat{\gamma}_n)$ be the value that maximizes
%
%
\begin{eqnarray} \label{eqllk3}
l_n(\beta,\gamma) &=& n^{-1} \sum_{i=1}^n \Biggl[\int\sum_{j=1}^{q_n}
\gamma_j B_j(t-X_i'\beta) \,dN_i(t) \nonumber\\[-7pt]\\[-7pt]
&&\hspace*{36pt}{} - \int I(Y_i\geq t) \exp\Biggl\{\sum_{j=1}^{q_n} \gamma_j
B_j(t-X_i'\beta)\Biggr\} \,dt\Biggr].\nonumber
\end{eqnarray}
Taking the first order derivatives of $l_n(\beta,\gamma)$ with respect
to $\beta$ and $\gamma$ and setting them to zero, we can obtain the
score equations. Since the integrals here are univariate integrals,
their numerical implementation can be easily done by the
one-dimensional Gaussian-quadrature method. A Newton--Raphson algorithm
or any other gradient-based search algorithm can be applied to solve
the score equations for all parameters $\theta= (\beta, \gamma)$, for example,
\[
\theta^{(m+1)} = \theta^{(m)} - H\bigl(\theta^{(m)}\bigr)^{-1}
\cdot S\bigl(\theta^{(m)}\bigr),
\]
where $\theta^{(m)}=(\beta^{(m)},\gamma^{(m)})$ is the parameter
estimate from the $m$th iteration, and
\[
S(\theta)= \pmatrix{
\displaystyle \frac{\partial l_n(\beta,\gamma)}{\partial\beta} \vspace*{2pt}\cr
\displaystyle \frac{\partial l_n(\beta,\gamma)}{\partial\gamma}},\qquad
H(\theta)= \pmatrix{
\displaystyle \frac{\partial^2 l_n(\beta,\gamma)}{\partial\beta\,\partial\beta'} &
\displaystyle \frac{\partial^2 l_n(\beta,\gamma)}{\partial\beta\,\partial\gamma'} \vspace*{2pt}\cr
\displaystyle \frac{\partial^2 l_n(\beta,\gamma)}{\partial\gamma\,\partial\beta'} &
\displaystyle \frac{\partial^2 l_n(\beta,\gamma)}{\partial\gamma\,\partial\gamma'}}
\]
are the score function and Hessian matrix of parameter $\theta$. For
any fixed $\beta$ and $n$, it is clearly seen that $l_n(\beta, \gamma)$
in (\ref{eqllk3}) is concave with respect to $\gamma$ and goes to
$-\infty$ if any $\gamma_j$ approaches either $\infty$ or $-\infty$;
hence $\hat\gamma_n$ must be bounded which yields an estimator of $s$
in $\mathcal{S}_n(T_{K_n},K_n,p)$.

As\vspace*{1pt} stated in the next section, the distribution of $\hat{\beta}_n$ can
be approximated by a normal distribution. One way to estimate the
variance matrix of $\hat{\beta}_n$ is to approximate the (inverse of
the) information matrix based on the efficient score function for $\beta
_0$ by plugging in the estimated parameters~$(\hat{\beta}_n, \hat
{\lambda}_n(\cdot))$. The consistency of such a variance estimator is
given in Theorem~\ref{thmvariance}. Another way is to invert the
observed information matrix from the last Newton--Raphson iteration,
taking into account that we are also estimating the nuisance parameter
$\gamma$. The consistency of the latter approach may be proved in a
similar way as Example 4 in \cite{Shen1997} or via Theorem~2.2 in \cite
{HeShao2000}, and we leave detailed derivation to interested readers.
Simulations indicate that both estimators work reasonably well.

\section{Asymptotic results} \label{secasymp}
Denote $\epsilon_{\beta}=Y-X'\beta$ and $\epsilon_0=Y-X'\beta_0$. We
assume the following regularity conditions:
\begin{longlist}[(C.2)]
\item[(C.1)] The true parameter $\beta_0$ belongs to the interior of a
compact set $\mathcal{B}\subseteq\mathbb{R}^d$.\vadjust{\goodbreak}

\item[(C.2)] (a) The covariate $X$ takes values in a bounded subset
$\mathcal{X} \subseteq\mathbb{R}^d$;\\ \hphantom{(C.2)\ \ \
\hspace*{9.2pt}}(b) $E(XX')$ is nonsingular.

\item[(C.3)] There is a truncation time $\tau<\infty$ such that, for
some constant $\delta$, $P(\epsilon_0 > \tau|X)\geq\delta>0$ almost
surely with respect to the probability measure of $X$. This implies
that $\Lambda_0(\tau) \leq-\log\delta<\infty$.

\item[(C.4)] The error $e_0$'s density $f$ and its derivative $\dot{f}$
are bounded and
\[
\int\bigl(\dot{f}(t)/f(t)\bigr)^2f(t) \,dt < \infty.
\]

\item[(C.5)] The conditional density of $C$ given $X$ and its
derivative $\dot{g}_{C|X}$ are uniformly bounded for all possible
values of $X$, that is,
\[
\sup_{x \in\mathcal{X}}g_{C|X}(t|X=x)\leq K_1,\qquad \sup_{x \in
\mathcal{X}}\bigl|\dot{g}_{C|X}(t|X=x)\bigr|\leq K_2
\]
for all $t\leq\tau$ with some constants $K_1, K_2>0$, where $\tau$ is
the truncation time defined in condition (C.3).

\item[(C.6)] Let $\mathcal{G}^p$ denote the collection of bounded
functions $g$ on $[a,b]$ with bounded derivatives $g^{(j)}$, $j=1,
\ldots
,k$, and the $k$th derivative $g^{(k)}$ satisfies the following
Lipschitz continuity condition:
\[
\bigl|g^{(k)}(s)-g^{(k)}(t)\bigr| \leq L|s-t|^m \qquad
\mbox{for } s, t \in[a,b],
\]
where $k$ is a positive integer and $m \in(0,1]$ such that $p=k+m \geq
3$, and $L<\infty$ is an unknown constant.
The true log hazard function $g_0(\cdot)=\log\lambda_0(\cdot)$ belongs
to~$\mathcal{G}^p$, where $[a,b]$ is a bounded interval.

\item[(C.7)] For some $\eta\in(0,1)$, $u'\operatorname{Var}(X|\epsilon_0)u \geq\eta
u'E(XX'|\epsilon_0)u$ almost surely for all $u\in\mathbb{R}^d$.
\end{longlist}

Condition (C.1) is a common regularity assumption that has been imposed
in the literature; see, for example, \cite{LaiYing1991}. Conditions
(C.2)(a), (C.3) and (C.4) were also assumed in \cite{Tsiatis}.
Condition (C.5) implies Condition $B$ in \cite{Tsiatis}.
In condition~(C.6), we require $p\geq3$ to provide desirable
controls of the spline approximation error rates of the first and
second derivatives of $g_0$ (see Corollary 6.21 of \cite{Schumaker}),
which are needed in verifying assumptions \mbox{(A4)--(A6)}.
Condition (C.7) was also proposed for the panel count data model in~\cite
{WellnerZhang2007}. As noted in their Remark 3.4, this condition~(C.7)
can be justified in many applications when condition (C.2)(b) is
satisfied. The bounded interval $[a,b]$ in (C.6) may be chosen as $a =
\inf_{y,x}(y-x'\beta_0) > -\infty$ and $b=\tau< \infty$ un\-der~\mbox{(C.1)--(C.3)},
which is what we use in the following.

Now define the collection of functions $\mathcal{H}^p$ as follows:
\[
\mathcal{H}^p = \{\zeta(\cdot,\beta)\dvtx\zeta(t,x,\beta) = g(\psi
(t,x,\beta)), g \in\mathcal{G}^p, t \in[a,b], x \in\mathcal{X},
\beta\in\mathcal{B} \},
\]
where
\[
\psi(t,x,\beta) = t-x'(\beta-\beta_0)
\]
and $\mathcal{G}^p$ is defined in (C.6). Here $\zeta$ is a composite
function of $g$ composed with~$\psi$.
Note that $\zeta(t,x,\beta_0) = g(t)$. Then for\vadjust{\goodbreak} $\zeta(\cdot,\beta) \in
\mathcal{H}^p$ we define the following norm:
%
%
\begin{equation} \label{eqnorm}
\| \zeta(\cdot,\beta) \|_2=\biggl\{\int_{\mathcal{X}}\int_{a}^{b} \bigl\{
g\bigl(t-x'(\beta-\beta_0)\bigr)\bigr\}^2 \,d \Lambda_0(t) \,d F_X(x)
\biggr\}^{1/2}.
\end{equation}
We also have the following collection of scores:
\[
\mathbb{H}=\biggl\{h\dvtx h(\cdot,\beta)=\frac{\partial\zeta_\eta(\cdot
,\beta)}{\partial\eta}\bigg|_{\eta=0} =w(\psi(\cdot,\beta)), \zeta
_\eta\in\mathcal{H}^p\biggr\},
\]
in which $h(t,x,\beta) = w(\psi(t,x,\beta))=w(t-x'(\beta-\beta_0))$.

For any $\theta_1=(\beta_1,\zeta_1(\cdot,\beta_1))$ and $\theta_2=(\beta
_2,\zeta_2(\cdot,\beta_2))$ in the space of $\Theta^p = \mathcal
{B}\times\mathcal{H}^p$, define the following distance:
%
%
\begin{equation} \label{eqdistance}
d(\theta_1,\theta_2)=\{|\beta_1-\beta_2|^2 + \|\zeta_1(\cdot,\beta
_1)-\zeta_2(\cdot,\beta_2)\|_2^2 \}^{1/2}.
\end{equation}
Let $\mathcal{G}^p_n =\mathcal{S}_n(T_{K_n}, K_n, p)$. Denote
\[
\mathcal{H}^p_n=\{\zeta(\cdot,\beta)\dvtx\zeta(t,x,\beta) = g(\psi
(t,x,\beta)), g \in\mathcal{G}^p_n, t \in[a,b], x \in\mathcal{X},
\beta\in\mathcal{B} \}
\]
and $\Theta^p_n = \mathcal{B}\times\mathcal{H}^p_n$. Clearly $\mathcal
{H}^p_n \subseteq\mathcal{H}_{n+1}^p \subseteq\cdots\subseteq
\mathcal{H}^p $ for all $n \geq1$. The sieve estimator $\hat{\theta
}_n=(\hat{\beta}_n,\hat{\zeta}_n(\cdot,\hat{\beta}_n))$, where $\hat
{\zeta}_n(t,x,\hat{\beta}_n) = \hat{g}_n(t-x'(\hat{\beta}_n-\beta_0))$,
is the maximizer of the empirical log-likelihood $n^{-1}l_n(\theta;Z)$
over the sieve space~$\Theta^p_n$. The following theorem gives the
convergence rate of the proposed estimator~$\hat{\theta}_n$ to the true
parameter $\theta_0=(\beta_0,\zeta_0(\cdot,\beta_0)) = (\beta_0, g_0)$.
%
%
%
\begin{theorem} \label{thmConvRate}
Let $K_n=O(n^{\nu})$, where $\nu$ satisfies the restriction $\frac
{1}{2(1+p)}<\nu<\frac{1}{2p}$ with $p$ being the smoothness parameter
defined in condition \textup{(C.6)}. Suppose conditions \textup{(C.1)--(C.7)}
hold, and the
failure time $T$ follows mo\-del~(\ref{eqmodel1}). Then
\[
d(\hat{\theta}_n,\theta_0) = O_p\bigl\{n^{-\min(p\nu, (1-\nu)/2)}\bigr\},
\]
where $d(\cdot,\cdot)$ is defined in (\ref{eqdistance}).
\end{theorem}
\begin{Remark*}
It is worth pointing out that the sieve space $\mathcal
{G}^p_n$ does not have to be restricted to the B-spline space; it can
be any sieve space as long as the estimator $\hat{\theta}_n \in\mathcal
{B}\times\mathcal{H}^p_n$ satisfies the conditions of Theorem 1 in
\cite{ShenWong1994}. We refer to \cite{chenxiaotong} for a
comprehensive discussion of the sieve estimation for semiparametric
models in general sieve spaces. Our choice of the B-spline space is
primarily motivated by its simplicity of numerical implementation,
which is a tremendous advantage of the proposed approach over exiting
numerical methods for the accelerated failure time models, in
particular, the linear programming approach.
\end{Remark*}

We provide a proof of Theorem \ref{thmConvRate} in the
supplementary material \cite{DingNan11} by checking the conditions of Theorem 1 in \cite
{ShenWong1994}. Theorem \ref{thmConvRate} implies that if $\nu
=1/(1+2p)$, $d(\hat{\theta}_n, \theta_0) = O_p(n^{-p/(1+2p)})$ which is
the optimal convergence rate in the nonparametric regression setting.
Although the overall convergence\vadjust{\goodbreak} rate is slower than $n^{-1/2}$, the
next theorem states that the proposed estimator of the regression
parameter is still asymptotically normal and semiparametrically efficient.
%
%
\begin{theorem} \label{thmnormality}
Given the following efficient score function for the censored linear
model derived by \cite{RitovWellner1988}:
\[
l_{\beta_0}^*(Y,\Delta,X)=\int\{X-P(X|Y-X'\beta_0 \geq t)\}\biggl\{
-\frac{\dot{\lambda}_0}{\lambda_0}(t)\biggr\}\, dM(t),
\]
where
\[
M(t)=\Delta I(Y-X'\beta_0 \leq t)-\int_{-\infty}^t I(Y-X'\beta_0 \geq
s)\lambda_0(s) \,ds
\]
is the failure counting process martingale, and
\[
P(X|Y-X'\beta_0 \geq t) = \frac{P\{XI(Y-X'\beta_0 \geq t)\}}{P\{
I(Y-X'\beta_0 \geq t)\}}
\]
was shown by \cite{Ritov}. Suppose that the conditions in Theorem \ref
{thmConvRate} hold, and $I(\beta_0)=P\{l_{\beta_0}^* (Y,\Delta
,X)^{\otimes2}\}$ is nonsingular, then
\[
n^{1/2}(\hat{\beta}_n-\beta_0)=n^{-1/2}I^{-1}(\beta_0)\sum_{i=1}^n
l_{\beta_0}^* (Y_i,\Delta_i,X_i) + o_p(1)
\rightarrow N(0,I^{-1}(\beta_0))
\]
in distribution.
\end{theorem}

The proof of Theorem \ref{thmnormality} is where we need to apply our
general sieve M-theorem proposed in Section \ref{secgenThm}. We prove
by checking assumptions \mbox{(A1)--(A6)}. Details are provided in Section \ref
{secProofs}. The following theorem gives consistency of the variance
estimator based on the above efficient score.
%
%
\begin{theorem} \label{thmvariance}
Suppose the conditions in Theorem \ref{thmnormality} hold. Denote
\[
l_{\hat\beta_n}^* (Y, \Delta, X) = \int\{X- \bar X(t; \hat\beta_n)\}
\{- \dot{\hat{g}}_n(t)\} \,d \hat{M}(t),
\]
where
\[
\bar X(t; \hat\beta_n) = \frac{\mathbb{P}_n \{ X I(Y - X'\hat\beta_n
\geq t)\}}{\mathbb{P}_n \{I(Y - X'\hat\beta_n \geq t)\}}
\]
and
\[
\hat{M}(t) = \Delta I(Y-X'\hat\beta_n \leq t)-\int_{-\infty}^t
I(Y-X'\hat\beta_n \geq s) \exp\{\hat{g}_n(s)\} \,ds.
\]
Then $\mathbb{P}_n \{l_{\hat\beta_n}^* (Y, \Delta, X)^{\otimes2} \}
\to P \{l^*_{\beta_0}(Y,\Delta,X)^{\otimes2} \} = I(\beta_0)$ in probability.
\end{theorem}

It is clearly seen that $\bar X(t, \hat\beta_n)$ in Theorem \ref
{thmvariance} estimates \mbox{$P(X|Y\,{-}\,X'\beta_0\,{\geq}\,t)$} in Theorem \ref
{thmnormality}. The proof of Theorem \ref{thmvariance} is provided in
the supplementary material~\cite{DingNan11}.\vadjust{\goodbreak}

\section{Numerical examples}
\subsection{Simulations} \label{secsimu}

Extensive simulations are carried out to evaluate the finite sample
performance of the proposed method. In the simulation studies, failure
times are generated from the model
\[
\log T = 2+X_1+X_2+e_0,
\]
where $X_1$ is Bernoulli with success probability 0.5, $X_2$ is
independent normal with mean 0 and standard deviation 0.5 truncated at
$\pm2$. This is the same model used by \cite{Jin2006} and \cite
{ZengLin2007}. We consider six error distributions: standard normal;
standard extreme-value; mixtures of $N(0,1)$ and $N(0,3^2)$ with mixing
probabilities $(0.5,0.5)$ and $(0.95,0.05)$, denoted by
$0.5N(0,1)+0.5N(0,3^2)$ and $0.95N(0,1)+0.05N(0,3^2)$, respectively;
Gumbel$(-0.5\mu,0.5)$ with $\mu$ being the Euler constant and
$0.5N(0,1)+0.5N(-1,0.5^2)$. The first four distributions were also
considered by \cite{ZengLin2007}.
Similar to \cite{ZengLin2007}, the censoring times are generated from
uniform $[0,c]$ distribution, where $c$ is chosen to produce a 25\%
censoring rate. We set the sample size $n$ to 200, 400 and 600.

%
\begin{table}[t!]
\tabcolsep=0pt
\caption{Summary statistics for the simulation studies. The true slope
parameters are $\beta_1=1$ and $\beta_2=1$. \textup{(a)}:~$N(0,1)$;
\textup{(b)}: standard extreme-value; \textup{(c)}: $0.5N(0,1)+0.5N(0,3^2)$;
\textup{(d)}:~$0.95N(0,1)+0.05N(0,3^2)$; \textup{(e)}: Gumbel($-0.5\mu$,0.5);
\textup{(f)}: $0.5N(0,1)+0.5N(-1,0.5^2)$}
\label{tabletab2}
{\fontsize{8.6pt}{11pt}\selectfont{
\begin{tabular*}{\tablewidth}{@{\extracolsep{\fill}}lccd{2.3}cccd{2.3}cd{2.3}cc@{}}
\hline
\multirow{2}{15pt}[-7pt]{\textbf{Err. dist}} & & &
\multicolumn{4}{c}{\textbf{B-spline MLE}}
& \multicolumn{2}{c}{\textbf{Log-rank}} & \multicolumn{2}{c}{\textbf{Gehan}} &
\\[-4pt]
& & & \multicolumn{4}{c}{\hspace*{-1.5pt}\hrulefill}
& \multicolumn{2}{c}{\hspace*{-1.5pt}\hrulefill}
& \multicolumn{2}{c}{\hspace*{-1.5pt}\hrulefill} & \\
& \multicolumn{1}{c}{$\bolds{n}$} & & \multicolumn{1}{c}{\textbf{Bias}}
& \multicolumn{1}{c}{\textbf{SE}} & \multicolumn{1}{c}{$\bolds{^1}$\textbf{SEE (CP)}}
& \multicolumn{1}{c}{$\bolds{^2}$\textbf{SEE (CP)}} & \multicolumn{1}{c}{\textbf{Bias}}
& \multicolumn{1}{c}{\textbf{SE}} & \multicolumn{1}{c}{\textbf{Bias}}
& \multicolumn{1}{c}{\textbf{SE}} & \multicolumn{1}{c@{}}{$\bolds{\sigma^*}$} \\
\hline
(a) & 200 & $\beta_1$ & 0.003 & 0.168 & 0.149 (0.912) & 0.155 (0.924) & 0.000
& 0.170 & 0.002 & 0.159 & 0.155 \\
& & $\beta_2$ & 0.003 & 0.167 & 0.153 (0.928) & 0.156 (0.928) & 0.004 & 0.171
& 0.002 & 0.160 & 0.156 \\
& 400 & $\beta_1$ & 0.006 & 0.110 & 0.108 (0.948) & 0.110 (0.950) & 0.005 &
0.115 & 0.008 & 0.108 & 0.110 \\
& & $\beta_2$ & 0.001 & 0.110 & 0.109 (0.944) & 0.110 (0.945) & 0.002 & 0.116
& 0.001 & 0.109 & 0.110 \\
& 600 & $\beta_1$ & 0.001 & 0.092 & 0.088 (0.939) & 0.090 (0.943) & 0.001 &
0.096 & 0.002 & 0.093 & 0.090 \\
& & $\beta_2$ & 0.005 & 0.091 & 0.089 (0.945) & 0.090 (0.944) & 0.005 & 0.097
& 0.003 & 0.092 & 0.090 \\[6pt]
(b) & 200 & $\beta_1$ & -0.009 & 0.180 & 0.154 (0.894) & 0.161 (0.903) &
-0.008 & 0.168 & -0.007 & 0.190 & 0.165 \\
& & $\beta_2$ & 0.004 & 0.182 & 0.162 (0.903) & 0.163 (0.915) & 0.005 & 0.170
& 0.005 & 0.195 & 0.169 \\
& 400 & $\beta_1$ & 0.000 & 0.126 & 0.113 (0.914) & 0.115 (0.923) & -0.001 &
0.124 & 0.000 & 0.143 & 0.117 \\
& & $\beta_2$ & 0.008 & 0.118 & 0.116 (0.934) & 0.116 (0.938) & 0.010 & 0.116
& 0.012 & 0.135 & 0.120 \\
& 600 & $\beta_1$ & 0.001 & 0.102 & 0.093 (0.919) & 0.094 (0.923) & 0.001 &
0.100 & 0.000 & 0.114 & 0.095 \\
& & $\beta_2$ & 0.011 & 0.098 & 0.095 (0.944) & 0.095 (0.945) & 0.011 & 0.097
& 0.007 & 0.114 & 0.098 \\[6pt]
(c ) & 200 & $\beta_1$ & 0.014 & 0.300 & 0.281 (0.930) & 0.279 (0.924) &
-0.020 & 0.315 & -0.019 & 0.292 & 0.259 \\
& & $\beta_2$ & 0.000 & 0.306 & 0.285 (0.916) & 0.282 (0.918) & 0.002 & 0.317
& 0.002 & 0.288 & 0.260 \\
& 400 & $\beta_1$ & 0.034 & 0.199 & 0.206 (0.955) & 0.200 (0.949) & 0.002 &
0.218 & 0.002 & 0.197 & 0.183 \\
& & $\beta_2$ & -0.003 & 0.207 & 0.208 (0.949) & 0.202 (0.942) & -0.001 &
0.222 & -0.002 & 0.200 & 0.184 \\
& 600 & $\beta_1$ & 0.035 & 0.168 & 0.171 (0.957) & 0.165 (0.949) & 0.003 &
0.185 & 0.001 & 0.163 & 0.150 \\
& & $\beta_2$ & -0.007 & 0.169 & 0.172 (0.956) & 0.166 (0.956) & -0.004 &
0.190 & -0.002 & 0.168 & 0.150 \\[6pt]
(d) & 200 & $\beta_1$ & -0.013 & 0.172 & 0.157 (0.926) & 0.164 (0.927) &
-0.010 & 0.181 & -0.007 & 0.166 & 0.167 \\
& & $\beta_2$ & -0.004 & 0.180 & 0.160 (0.908) & 0.164 (0.913) & -0.005 &
0.184 & -0.005 & 0.173 & 0.166 \\
& 400 & $\beta_1$ & 0.003 & 0.119 & 0.113 (0.944) & 0.116 (0.948) & 0.004 &
0.126 & 0.006 & 0.117 & 0.118 \\
& & $\beta_2$ & 0.003 & 0.117 & 0.114 (0.942) & 0.116 (0.953) & 0.004 & 0.126
& 0.003 & 0.115 & 0.118 \\
& 600 & $\beta_1$ & -0.003 & 0.097 & 0.093 (0.948) & 0.095 (0.952) & -0.002 &
0.105 & 0.002 & 0.097 & 0.096 \\
& & $\beta_2$ & 0.001 & 0.096 & 0.094 (0.942) & 0.095 (0.944) & 0.002 & 0.105
& 0.003 & 0.094 & 0.096 \\[6pt]
(e) & 200 & $\beta_1$ & 0.004 & 0.080 & 0.077 (0.944) & 0.078 (0.946) &
-0.001 & 0.111 & 0.004 & 0.088 & 0.079 \\
& & $\beta_2$ & -0.001 & 0.083 & 0.080 (0.929) & 0.078 (0.934) & 0.000 & 0.114
& 0.000 & 0.091 & 0.080 \\
& 400 & $\beta_1$ & -0.005 & 0.055 & 0.055 (0.946) & 0.055 (0.951) & -0.003 &
0.079 & -0.004 & 0.061 & 0.056 \\
& & $\beta_2$ & 0.003 & 0.055 & 0.056 (0.954) & 0.056 (0.950) & 0.003 & 0.081
& 0.003 & 0.063 & 0.056 \\
& 600 & $\beta_1$ & -0.003 & 0.047 & 0.045 (0.940) & 0.045 (0.938) & 0.000 &
0.067 & -0.001 & 0.052 & 0.045 \\
& & $\beta_2$ & -0.001 & 0.047 & 0.046 (0.944) & 0.045 (0.943) & -0.002 &
0.066 & -0.001 & 0.051 & 0.046 \\[6pt]
(f) & 200 & $\beta_1$ & -0.002 & 0.126 & 0.117 (0.918) & 0.120 (0.929) &
-0.002 & 0.159 & -0.001 & 0.128 & 0.119 \\
& & $\beta_2$ & 0.000 & 0.133 & 0.120 (0.917) & 0.121 (0.926) & 0.002 & 0.164
& 0.001 & 0.134 & 0.116 \\
& 400 & $\beta_1$ & -0.002 & 0.087 & 0.084 (0.949) & 0.085 (0.950) & 0.003 &
0.114 & 0.000 & 0.091 & 0.084 \\
& & $\beta_2$ & 0.004 & 0.086 & 0.086 (0.951) & 0.086 (0.953) & 0.003 & 0.111
& 0.004 & 0.090 & 0.082 \\
& 600 & $\beta_1$ & 0.003 & 0.074 & 0.070 (0.929) & 0.070 (0.931) & 0.005 &
0.101 & 0.001 & 0.074 & 0.069 \\
& & $\beta_2$ & 0.003 & 0.074 & 0.070 (0.936) & 0.070 (0.936) & 0.009 & 0.104
& 0.004 & 0.075 & 0.067 \\
\hline
\end{tabular*}}}
\end{table}

We choose cubic B-splines with one interior knot for $n=200$ and $400$,
and two interior knots for $n=600$. We perform the sieve maximum
likelihood analysis and obtain the estimates of the slope parameters
using the Newton--Raphson algorithm that updates $(\beta, \gamma)$
iteratively. We stop iteration when the change of parameter estimates
or the gradient value is less than a pre-specified tolerance value that
is set to be $10^{-5}$ in our simulations. Log-rank and Gehan-weighted
estimators are included for efficiency comparisons. We calculate the
theoretical semiparametric efficiency bound $I^{-1}(\beta_0)$, and
scale it by the sample size, that is, $\sigma^*=\sqrt{I^{-1}(\beta
_0)/n}$, which serves as the reference standard error under the fully
efficient situation. Table \ref{tabletab2} summarizes the results of
these studies based on 1,000 simulated datasets. The bias of the
proposed estimators of $\beta_1$ and $\beta_2$ are negligible. Both
variance estimation procedures, denoted as $^1$SEE (the standard error
estimates by inverting the information matrix based on the efficient
score function) and $^2$SEE (the standard error estimates by inverting
the observed information matrix of all parameters including nuisance
parameters), yield nice standard error estimates for the parameter
estimators comparing to the empirical standard error SE, and the $95\%$
confidence intervals have proper coverage probabilities, especially
when the sample size is large. For the $N(0,1)$ error and the two
mixtures of normal errors that are also considered in \cite
{ZengLin2007}, the proposed estimators are more efficient than the
log-rank estimators and have similar variances to the Gehan-weighted
estimators. For the standard extreme-value error, the proposed
estimators are more efficient than the Gehan-weighted estimator and
similar to the log-rank estimator that is known to be the most
efficient estimator under this particular error distribution. For the
Gumbel$(-0.5\mu,0.5)$ and $0.5N(0,1)+0.5N(-1,0.5^2)$ errors, the
proposed estimators
are more efficient than the other two estimators. Under all six
error distributions, the standard errors of the proposed estimators are
close to the efficient theoretical standard errors.
The sample averages of the estimates for $\lambda_0$ under different
simulation settings are reasonably close to corresponding true curves
(results not shown here; see \cite{Ding} for details).

\subsection{A real data example} \label{secexample}
We use the Stanford heart transplant data \cite{Miller1982} as an
illustrative example. This dataset was also analyzed by \cite{Jin2006}
using their proposed least squares estimators. Following their
analysis, we consider the same two models: the first one regresses the
base-10 logarithm of the survival time on age at transplant and T5
mismatch score for the 157 patients with complete records on T5
measure, and the second one regresses the base-10 logarithm of the
survival time on age and age$^2$. There were 55 censored patients.
We fit these two models using the proposed method with five cubic
B-spline basis functions.

%
\begin{table}
\caption{Regression parameter estimates and standard error estimates
for the Stanford heart transplant data. The proposed estimators are
compared with Gehan-weighted estimators reported in \protect\cite{Jin2003} and
Buckley--James estimators reported in \protect\cite{Miller1982}}
\label{tabletab4}
%
%
\begin{tabular*}{\tablewidth}{@{\extracolsep{\fill}}lcd{2.4}cd{2.4}cd{2.4}d{1.4}@{}}
\hline
& & \multicolumn{2}{c}{\textbf{B-spline MLE}} &
\multicolumn{2}{c}{\textbf{Gehan-weighted}} &
\multicolumn{2}{c@{}}{\textbf{Buckley--James}}
\\[-4pt]
& & \multicolumn{2}{c}{\hspace*{-1.5pt}\hrulefill} &
\multicolumn{2}{c}{\hspace*{-1.5pt}\hrulefill} &
\multicolumn{2}{c@{}}{\hspace*{-1.5pt}\hrulefill}
\\
& \multicolumn{1}{c}{\textbf{Covariate}} & \multicolumn{1}{c}{\textbf{Est.}}
& \multicolumn{1}{c}{\textbf{SE}} & \multicolumn{1}{c}{\textbf{Est.}}
& \multicolumn{1}{c}{\textbf{SE}} & \multicolumn{1}{c}{\textbf{Est.}}
& \multicolumn{1}{c@{}}{\textbf{SE}} \\
\hline
M.1 & Age & -0.0237 & 0.0068 & -0.0211 & 0.0106 & -0.015 & 0.008
\\
& T5 & -0.2118 & 0.1271 & -0.0265 & 0.1507 & -0.003 & 0.134
\\
M.2 & Age & 0.1022 & 0.0245 & 0.1046 & 0.0474 & 0.107 & 0.037 \\
& Age$^2$ & -0.0016 & 0.0004 & -0.0017 & 0.0006 & -0.0017 & 0.0005
\\
\hline
\end{tabular*}
\end{table}

We report the parameter estimates and the standard error estimates in
Table~\ref{tabletab4} and compare them with the Gehan-weighted
estimators reported by \cite{Jin2006} and the Buckley--James estimators
reported by \cite{Miller1982}. For the first model, the parameter
estimates for the age effect are fairly similar among all estimators,
and the standard error estimate from the proposed method tends to be
smaller, while the parameter estimates for the T5 mismatch score vary
across different estimators with none of them being significant at the
0.05 level. The disparity of the T5 effect may be due to what was
pointed out by \cite{Miller1982}: the accelerated failure time model
with age and T5 as covariates does not fit the data ideally. For the
second model with age and age$^2$ being the covariates, the point
estimates are very similar across all methods and the standard error
estimates from the proposed method are the smallest.

\section{Discussion} \label{secdisc}
By applying the proposed general sieve M-estimation theory for
semiparametric models with bundled parameters, we are able to derive
the asymptotic distribution for the sieve maximum likelihood estimator
in a linear regression model where the response variable is subject to
right censoring.
By providing a both statistically and computationally efficient
estimating procedure, this work makes the linear model a more viable
alternative to the Cox proportional hazards model. Comparing to the
existing methods for estimating $\beta$ in a linear model, the proposed
method has three advantages. First, the estimating functions are smooth
functions in contrast to the discrete estimating functions in the
existing estimation methods; thus the root search is easier and can be
done quickly by conventional iterative methods such as the
Newton--Raphson algorithm. Second, the standard error estimates are
obtained directly by inverting either the efficient information matrix
for the regression parameters or the observed information matrix of all
parameters; either method is more computationally tractable compared to
the re-sampling techniques. Third, the proposed estimator achieves the
semiparametric efficiency bound.

The proposed general sieve M-estimation theory can also be applied to
other statistical models, for example, the single index model, the Cox
model with an unknown link function and the linear model under
different censoring mechanisms. Such research is undergoing and will be
presented elsewhere.

\section{\texorpdfstring{Proof of Theorem \protect\ref{thmnormality}}{Proof of Theorem 4.2}}\label{secProofs}
Empirical process theory developed in \cite{Vaart1998,vanWellner1996}
will be heavily involved in the proof. We use the symbol $\lesssim$ to
denote that the left-hand side is bounded above by a constant times the
right-hand side and~$\gtrsim$ to denote that the left-hand side is
bounded below by a constant times the right-hand side. For notational
simplicity, we drop the superscript $*$ in the outer probability
measure $P^*$ whenever an outer probability applies.

\subsection{Technical lemmas}
We first introduce several lemmas that will be used for the proofs of
Theorems \ref{thmConvRate}, \ref{thmnormality} and \ref
{thmvariance}. Proofs of these lemmas are provided in the
supplementary material \cite{DingNan11}.
%
%
\begin{lem} \label{lemlem1}
Under conditions \textup{(C.1)--(C.3)} and \textup{(C.6)}, the log-likelihood
\begin{eqnarray*}
l(\beta,\zeta(\cdot,\beta);Z)
&=& \Delta g\bigl(\epsilon_0-X'(\beta-\beta_0)\bigr)\\
&&{}-\int_a^b 1(\epsilon_0
\geq t)\exp\bigl\{g\bigl(t-X'(\beta-\beta_0)\bigr)\bigr\} \,dt,
\end{eqnarray*}
where $\epsilon_0 = Y-X'\beta_0$, has bounded and continuous first and
second derivatives with respect to $\beta\in\mathcal{B}$ and $\zeta
(\cdot,\beta) \in\mathcal{H}^p$.
\end{lem}
%
%
\begin{lem} \label{lemlem2}
For $g_0 \in\mathcal{G}^p$, there exists a function $g_{0,n} \in
\mathcal{G}^p_n$ such that
\[
\|g_{0,n} - g_0\|_{\infty} = O(n^{-p\nu}).\vadjust{\goodbreak}
\]
\end{lem}
%
%
\begin{lem} \label{lemlem3}
Let\vspace*{1pt} $\theta_{0,n}=(\beta_0,\zeta_{0,n}(\cdot,\beta_0))$ with $\zeta
_{0,n}(\cdot,\beta_0) \equiv g_{0,n}$ defined in Lemma \ref{lemlem2}.
Denote $\mathcal{F}_n = \{l(\theta;z)-l(\theta_{0,n};z)\dvtx\theta\in
\Theta^p_n\}$. Assume that conditions \textup{(C.1)--(C.3)} and
\textup{(C.6)} hold, then the $\varepsilon$-bracketing number
associated with \mbox{$\| \cdot\|_{\infty}$} norm for $\mathcal{F}_n$ is
bounded by $(1/\varepsilon)^{cq_n+d}$, that is, $N_{[\csp
]}(\varepsilon,\mathcal{F}_n,\|\cdot\|_{\infty})
\lesssim(1/\varepsilon)^{cq_n+d}$ for some constant $c>0$.
\end{lem}
%
%
\begin{lem} \label{lemlem4}
Let\vspace*{1pt} $h_j^*(t,x,\beta)=w_j^*(\psi(t,x,\beta))$, where
$h_j^*(t,x,\beta _0)=w_j^*(t) = -\dot{g}_0(t)P(X_j|\epsilon_0 \geq t),
j=1,\ldots,d$. Assume conditions \textup{(C.1)--(C.6)} hold, then there
exists $h_{j,n}^*(t,x,\beta) = w_{j,n}^*(\psi(t,x,\beta))
\in\mathcal{H}^2_n$ such that $\| h_{j,n}^* - h_j^*\|_{\infty} =
O(n^{-2\nu})$, or equivalently, $\| w_{j,n}^* - w_j^*\|_{\infty} =
O(n^{-2\nu})$ where $w_{j,n}^* \in \mathcal{G}_n^2$.
\end{lem}
%
%
\begin{lem} \label{lemlem5}
For $h_j^*$ defined in Lemma \ref{lemlem4}, denote the class of
functions
\[
\mathcal{F}_n^j(\eta)=\{\dot{l}_{\zeta}(\theta;z)[h_j^*-h_j]\dvtx\theta
\in
\Theta^p_n, h_j\in\mathcal{H}^2_n, d(\theta,\theta_0) \leq\eta, \|
h_j-h_j^*\|_{\infty} \leq\eta\}.
\]
Assume conditions \textup{(C.1)--(C.6)} hold, then
$N_{[\csp]}(\varepsilon,\mathcal{F}_n^j(\eta),\|\cdot\|_{\infty})
\lesssim(\eta/\varepsilon)^{cq_n+d}$ for some constant $c>0$.
\end{lem}
%
%
\begin{lem} \label{lemlem6}
For $j=1,\ldots,d$, define the following two classes of functions:
\begin{eqnarray*}
\mathcal{F}_{n,j}^{\beta}(\eta) & = & \{\dot{l}_{\beta
_j}(\theta;z)-\dot{l}_{\beta_j}(\theta_0;z)\dvtx
\theta\in\Theta^p_n, d(\theta,\theta_0) \leq\eta,
\\
&&\hspace*{45.8pt} \|\dot{g}(\psi(\cdot,\beta
))-\dot{g}_{0}(\psi(\cdot,\beta_0))\|_2 \leq\eta\}
\end{eqnarray*}
and
\[
\mathcal{F}_{n,j}^{\zeta}(\eta) = \{\dot{l}_{\zeta}(\theta
;z)[h_j^*(\cdot, \beta)]-\dot{l}_{\zeta}(\theta_0;z)[h_j^*(\cdot, \beta
_0)]\dvtx\theta\in\Theta^p_n, d(\theta,\theta_0) \leq\eta\},
\]
where $\dot{l}_{\beta_j}(\theta;Z)$ is the $j$th element of $\dot
{l}_{\beta}(\theta;Z)$, $\dot g(\cdot)$ denotes the derivative of
$g(\cdot)$ and $h_j^*$ is defined in Lemma \ref{lemlem5}. Assume
conditions \textup{(C.1)--(C.6)} hold, then $N_{[\csp]}(\varepsilon,\mathcal
{F}_{n,j}^{\beta}(\eta), \|\cdot\|_{\infty}) \lesssim(\eta/\varepsilon
)^{c_1q_n+d}$ and $N_{[\csp]}(\varepsilon,\mathcal{F}_{n,j}^{\zeta}(\eta
),\| \cdot\|_{\infty}) \lesssim(\eta/\varepsilon)^{c_2q_n+d}$ for
some constants $c_1,c_2>0$.
\end{lem}

\subsection{\texorpdfstring{Proof of Theorem \protect\ref{thmnormality}}{Proof of Theorem 4.2}}
We prove the theorem by checking assumptions (A1)--(A6) in Section \ref
{secgenThm}. Here the criterion function of a single observation is
the log-likelihood function $l(\beta,\zeta(\cdot,\beta);Z)$. So instead
of $m$, we use $l$ to denote the criterion function.
By Theorem \ref{thmConvRate} we know that assumption (A1) holds with $\xi
=\min(p\nu,(1-\nu)/2)$ and the norm $\| \cdot\|_2$ defined in (\ref
{eqnorm}). Assumption (A2) automatically holds for the scores. For (A3), we need to
find an $\mathbf{h}^* = (h_1^*,\ldots,h_d^*)'$ with $\mathbf
{h}^*(t,x,\beta_0)=\mathbf{w}^*(t)$ such that
\begin{eqnarray*}
&& \ddot{S}_{\beta\zeta}(\beta_0,\zeta_0(\cdot,\beta_0))[h]-\ddot
{S}_{\zeta\zeta}(\beta_0,\zeta_0(\cdot,\beta_0))[\mathbf{h}^*,h]\\
&&\qquad = P\{\ddot{l}_{\beta\zeta} (\beta_0,\zeta_0(\cdot,\beta
_0);Z)[h]-\ddot{l}_{\zeta\zeta}(\beta_0,\zeta_0(\cdot,\beta
_0);Z)[\mathbf{h}^*,h]\}=0
\end{eqnarray*}
for all $h \in\mathbb{H}$ with $h(t,x,\beta) = w(t-x'(\beta-\beta
_0))$. Note that
\begin{eqnarray*}
&& P\{\ddot{l}_{\beta\zeta}(\beta_0,\zeta_0(\cdot,\beta_0);Z)[h]-\ddot
{l}_{\zeta\zeta}(\beta_0,\zeta_0(\cdot,\beta_0);Z)[\mathbf{h}^*,h]\} \\
&&\qquad = P\biggl\{-X \biggl[\Delta\dot{w}(\epsilon_0)-\int
_{a}^{b}1(\epsilon_0 \geq t) \exp\{g_0(t)\}\dot{w}(t) \,dt \biggr] \\
&&\qquad\quad\hspace*{14pt}{} + \int_{a}^{b}1(\epsilon_0\geq t) \exp\{g_0(t)\}
w(t)[X\dot{g}_0(t)+\mathbf{w}^*(t)] \,dt \biggr\}.
\end{eqnarray*}
Since $P\{\dot{l}_{\zeta}(\beta_0,\zeta_0(\cdot,\beta_0);Z)[h]\vert X\}
= 0$ for all $h \in\mathbb{H}$, replacing $h(\cdot, \beta_0)$ by $\dot
{w}$ we have
\begin{eqnarray*}
&& P \biggl\{-X \biggl[\Delta\dot{w}(\epsilon_0)-\int_{a}^{b}1(\epsilon
_0 \geq t) \exp\{g_0(t)\}\dot{w}(t) \,dt \biggr]\biggr\} \\
&&\qquad = P\biggl\{-X \cdot P\biggl[\Delta\dot{w}(\epsilon_0)-\int
_{a}^{b}1(\epsilon_0 \geq t) \exp\{g_0(t)\}\dot{w}(t) \,dt \Big\vert X
\biggr] \biggr\}\\
&&\qquad = P\{-X \cdot0\} = 0.
\end{eqnarray*}
Hence we only need to find a $\mathbf{w}^*$ such that
\begin{eqnarray*}
&& P\biggl\{\int_{a}^{b}1(\epsilon_0 \geq t) \exp\{
g_0(t)\} w(t)[X\dot{g}_0(t)+\mathbf{w}^*(t)] \,dt \biggr\} \\
&&\qquad = \int_{a}^{b} \exp\{g_0(t)\}w(t) \{\dot
{g}_0(t)P[1(\epsilon_0 \geq t)X]+\mathbf{w}^*(t)P[1(\epsilon_0 \geq
t)]\} \,dt = 0.
\end{eqnarray*}
One obvious choice for $\mathbf{w}^*$ (or $\mathbf{h}^*$) is
%
%
\begin{equation} \label{eqhstar}\quad
\mathbf{h}^*(t,x,\beta_0) = \mathbf{w}^*(t)= -\dot{g}_0(t)\frac
{P[1(\epsilon_0 \geq t)X]}{P[1(\epsilon_0 \geq t)]} = -\dot
{g}_0(t)P(X|\epsilon_0 \geq t).
\end{equation}
Then it follows
\begin{eqnarray*}
&& \dot{l}_\beta(\beta_0,\zeta_0(\cdot,\beta_0);Z)-\dot{l}_{\zeta}(\beta
_0,\zeta_0(\cdot,\beta_0);Z)[\mathbf{h}^*] \\
&&\qquad = \Delta\{-\dot{g}_0(Y-X'\beta_0)\} \{X-P(X|\epsilon_0
\geq Y-X'\beta_0)\} \\
&&\qquad\quad{} - \int1(Y-X'\beta_0 \geq t) \{X-P(X|\epsilon_0
\geq t)\}\{-\dot{g}_0(t)\}\exp\{g_0(t)\} \,dt \\
&&\qquad = \int\{X-P(X|\epsilon_0 \geq t)\}\{-\dot{g}_0(t)\}\,
d M(t)\\
&&\qquad = l_{\beta_0}^*(Y,\Delta,X),
\end{eqnarray*}
which is the efficient score function for $\beta_0$ originally derived
by \cite{RitovWellner1988}, where
\[
M(t)=\Delta I(Y-X'\beta_0 \leq t)-\int_{-\infty}^t I(Y-X'\beta_0 \geq
s)\exp\{g_0(s)\} \,ds.
\]
By the fact of zero-mean for a score function, it is straightforward to
verify the following equalities:
\begin{eqnarray*}
P\ddot{l}_{\beta\zeta}(\beta,\zeta(\cdot,\beta);Z)[h] &=& - P\{
\dot{l}_\beta(\beta,\zeta(\cdot,\beta);Z)\dot{l}_\zeta'(\beta,\zeta
(\cdot,\beta);Z)[h]\}, \\
P\ddot{l}_{\zeta\beta}(\beta,\zeta(\cdot,\beta);Z)[h] &=& - P\{\dot
{l}_\zeta(\beta,\zeta(\cdot,\beta);Z)[h]\dot{l}_{\beta}'(\beta,\zeta
(\cdot,\beta);Z)\}, \\
P\ddot{l}_{\beta\beta}(\beta,\zeta(\cdot,\beta);Z)&=&-P\{\dot
{l}_\beta(\beta,\zeta(\cdot,\beta);Z)\dot{l}_\beta'(\beta,\zeta(\cdot
,\beta);Z)\}, \\
P\ddot{l}_{\zeta\zeta}(\beta,\zeta(\cdot,\beta);Z)[h_1,h_2]&=&-P\{
\dot{l}_\zeta(\beta,\zeta(\cdot,\beta);Z)[h_1] \dot{l}_\zeta'(\beta
,\zeta(\cdot,\beta);Z)[h_2]\}.
\end{eqnarray*}
Then together with the fact that
\[
P\{\ddot{l}_{\beta\zeta}(\beta_0,\zeta_0(\cdot,\beta_0);Z)[\mathbf
{h}^*]-\ddot{l}_{\zeta\zeta}(\beta_0,\zeta_0(\cdot,\beta_0);Z)[\mathbf
{h}^*,\mathbf{h}^*]\}=0,
\]
the matrix $A$ in assumption (A3) of Theorem \ref{thmgenThm} is given by
\begin{eqnarray*}
A&=& P\{-\ddot{l}_{\beta\beta}(\beta_0,\zeta_0(\cdot,\beta_0);Z)+\ddot
{l}_{\zeta\beta}(\beta_0,\zeta_0(\cdot,\beta_0);Z)[\mathbf{h}^*] \\
&&\hspace*{11.4pt}{} + \ddot{l}_{\beta\zeta}(\beta_0,\zeta_0(\cdot,\beta
_0);Z)[\mathbf{h}^*]-\ddot{l}_{\zeta\zeta}(\beta_0,\zeta_0(\cdot,\beta
_0);Z)[\mathbf{h}^*,\mathbf{h}^*]\}\\
&=& P\{\dot{l}_{\beta}(\beta_0,\zeta_0(\cdot,\beta_0);Z)\dot{l}_{\beta
}'(\beta_0,\zeta_0(\cdot,\beta_0);Z)\\
&&\hspace*{11.4pt}{} - \dot{l}_{\zeta}(\beta_0,\zeta_0(\cdot,\beta_0);Z)[\mathbf
{h}^*]\dot{l}_{\beta}'(\beta_0,\zeta_0(\cdot,\beta_0);Z) \\
&&\hspace*{11.4pt}{} - \dot{l}_{\beta}(\beta_0,\zeta_0(\cdot,\beta_0);Z)\dot
{l}_{\zeta}'(\beta_0,\zeta_0(\cdot,\beta_0);Z)[\mathbf{h}^*] \\
&&\hspace*{11.4pt}{}
+ \dot{l}_{\zeta}(\beta_0,\zeta_0(\cdot,\beta_0);Z)[\mathbf
{h}^*]\dot{l}_{\zeta}'(\beta_0,\zeta_0(\cdot,\beta_0);Z)[\mathbf{h}^*]\}
\\
&=& P\{\dot{l}_{\beta}(\beta_0,\zeta_0(\cdot,\beta_0);Z)-\dot{l}_{\zeta
}(\beta_0,\zeta_0(\cdot,\beta_0);Z)[\mathbf{h}^*]\}^{\otimes2}\\
&=&
Pl_{\beta_0}^*(Y,\Delta,X)^{\otimes2},
\end{eqnarray*}
which is the information matrix for $\beta_0$.

To verify (A4), we note that the first part automatically holds since
$\hat{\beta}_n$ satisfies the score equation $\dot{S}_{\beta,n}(\hat
{\beta}_n,\hat{\zeta}_n(\cdot, \hat\beta_n)) = \mathbb{P}_n \dot
{l}_{\beta}(\hat{\beta}_n,\hat{\zeta}_n(\cdot, \hat\beta_n);Z)=0$.
Next we shall show that
\begin{eqnarray*}
&& \dot{S}_{\zeta,n}(\hat{\beta}_n,\hat{\zeta}_n(\cdot
, \hat\beta_n))[h_j^*] \\
&&\qquad = \mathbb{P}_n \biggl\{\Delta w_j^*(Y-X'\hat{\beta
}_n)-\int1(Y\geq t) \exp\{\hat{\zeta}_n(t,X,\hat{\beta}_n)\}
w_j^*(t-X'\hat{\beta}_n) \,dt \biggr\} \\
&&\qquad = o_p(n^{-1/2}),
\end{eqnarray*}
where $w_j^*(t)=-\dot{g}_0(t)P(X_j|\epsilon_0 \geq t)$, $j=1,\ldots,d$,
is the $j$th component of~$\mathbf{w}^*(t)$ given in (\ref{eqhstar}).
According to Lemma \ref{lemlem4}, there exists $h_{j,n}^*\in\mathcal
{H}^2_n$ such that $\Vert h_j^*-h_{j,n}^* \Vert_{\infty} = O(n^{-2\nu
})$. Then by the score equation for $\gamma\dvtx\dot{S}_{\gamma,n}(\hat
{\beta}_n,\hat{\gamma}_n) = \mathbb{P}_n \dot{l}_{\gamma}(\hat{\beta
}_n,\hat{\gamma}_n;Z)=0$ and the fact that $w_{j,n}^*(t)$ can be
written as $w_{j,n}^*(t)=\sum_{k=1}^{q_n}\gamma_{j,k}^*B_k(t)$ for
some\vadjust{\goodbreak}
coefficients $\{\gamma_{j,1}^*,\ldots,\gamma_{j,q_n}^*\}$ and the\vspace*{1pt} basis
func-\break tions~$B_k(t)$ of the spline space, it follows that
\[
\mathbb{P}_n \biggl\{\Delta w_{j,n}^*(Y-X'\hat{\beta}_n)-\int1(Y\geq
t)\exp\{\hat{\zeta}_n(t,X,\hat{\beta}_n)\}w_{j,n}^*(t-X'\hat{\beta
}_n) \,dt\biggr\}=0.
\]
So it suffices to show that for each $1\leq j \leq d$,
\[
I_n = \mathbb{P}_n\dot{l}_\zeta(\hat{\beta}_n,\hat{\zeta}_n(\cdot,\hat
{\beta}_n);Z)[h_j^*-h_{j,n}^*] = o_p(n^{-1/2}).
\]
Since $P\{\dot{l}_\zeta(\beta_0,\zeta_0(\cdot,\beta
_0);Z)[h_j^*\,{-}\,h_{j,n}^*] \}\,{=}\,0$, we decompose $I_n$ into
\mbox{$I_n\,{=}\,I_{1n}\,{+}\,I_{2n}$}, where
\[
I_{1n} = (\mathbb{P}_n-P)\dot{l}_\zeta(\hat{\beta}_n,\hat{\zeta}_n(\cdot
,\hat{\beta}_n);Z)[h_j^*-h_{j,n}^*]
\]
and
\[
I_{2n} = P\{\dot{l}_\zeta(\hat{\beta}_n,\hat{\zeta}_n(\cdot,\hat
{\beta}_n);Z)[h_j^*-h_{j,n}^*] -\dot{l}_\zeta(\beta_0,\zeta_0(\cdot
,\beta_0);Z)[h_j^*-h_{j,n}^*] \}.
\]
We will show that $I_{1n}$ and $I_{2n}$ are both $o_p(n^{-1/2})$.

First consider $I_{1n}$. According to Lemma \ref{lemlem5}, the
$\varepsilon$-bracketing number associated with $\| \cdot\|_{\infty}$
norm for the class $\mathcal{F}_n^j(\eta)$ defined in Lemma \ref{lemlem5}
is bounded by $(\eta/\varepsilon)^{cq_n+d}$.
This implies that
\[
\log N_{[\csp]}(\varepsilon,\mathcal{F}_n^j(\eta),L_2(P)) \leq\log
N_{[\csp]}(\varepsilon,\mathcal{F}_n^j(\eta),\| \cdot\|_{\infty
})\lesssim q_n \log(\eta/\varepsilon),
\]
which leads to the bracketing integral
\begin{eqnarray*}
J_{[\csp]}(\eta,\mathcal{F}_n^j(\eta),L_2(P)) &=& \int_0^{\eta} \sqrt
{1+\log N_{[\csp]}(\varepsilon,\mathcal{F}_n^j(\eta),L_2(P))}\,
d\varepsilon\\
&\lesssim& q_n^{1/2}\eta.
\end{eqnarray*}
Now we pick $\eta$ to be $\eta_n = O\{n^{-\min(2\nu, (1-\nu)/2)}\}$, then
\[
\| h_j^*-h_{j,n}^* \|_{\infty} = O(n^{-2\nu}) \leq O\bigl\{n^{-\min(2\nu,
(1-\nu)/2)}\bigr\} = \eta_n,
\]
and since $p \geq3$,
\[
d(\hat{\theta}_n,\theta_0) = O_p\bigl\{n^{-\min(p\nu, (1-\nu)/2)}\bigr\} \leq
O_p\bigl\{n^{-\min(2\nu, (1-\nu)/2)}\bigr\} = \eta_n.
\]
Therefore, $\dot{l}_\zeta(\hat{\beta}_n,\hat{\zeta}_n(\cdot,\hat{\beta
}_n);z)[h_j^*-h_{j,n}^*] \in\mathcal{F}_n^j(\eta_n)$. Denote $t_{\beta
}=t-X'(\beta-\beta_0)$ for notational simplicity, for any $\dot{l}_\zeta
(\theta;Z)[h_j^*-h] \in\mathcal{F}_n^j(\eta_n)$, it follows that
\begin{eqnarray*}
&& P\{\dot{l}_\zeta(\theta;Z)[h_j^*-h]\}^2 \\
&&\qquad = P\biggl\{\Delta(w_j^*-w)(\epsilon_\beta) + \int_a^b
1(\epsilon_0 \geq t)\exp\{g(t_{\beta})\}(w_j^*-w)(t_{\beta}) \,dt
\biggr\}^2 \\
&&\qquad \lesssim\| w_j^*-w \|_{\infty}^2 + P\biggl\{\int_a^b \exp\{
2g(t_{\beta})\}(w_j^*-w)^2(t_{\beta}) \,dt\biggr\} \\
&&\qquad \lesssim\| w_j^*-w \|_{\infty}^2 + \| w_j^*-w \|_{\infty}^2
\int_a^b P[\exp\{2g(t_{\beta})\}] \,dt ,
\end{eqnarray*}
where the first inequality holds because of the Cauchy--Schwarz
inequality. Since $\|w_j^* - w\|_\infty\leq\eta_n$, by the same
argument as (\cite{ShenWong1994}, page 591), for slowly growing $c_n$
(their $l_n$), for example, $c_n=o(\log(\eta_n^{-1}))$, we know that
$\| \dot{l}_\zeta(\theta;Z)[h_j^*-h] \|_{\infty}$ is bounded by some
constant $0<M<\infty$ and $P\{\dot{l}_\zeta(\theta;Z)[h_j^*-h]\}^2
\lesssim\eta_n$ for a slightly enlarged $\eta_n$ obtained by a fine
adjustment of $\nu$. Then by the maximal inequality in Lemma 3.4.2 of
\cite{vanWellner1996}, it follows that
\begin{eqnarray*}
E_P\Vert\mathbb{G}_n \|_{\mathcal{F}_n^j(\eta_n)}
&\lesssim& J_{[\csp
]}(\eta_n,\mathcal{F}_n^j(\eta_n),L_2(P)) \biggl(1+\frac{J_{[\csp]}(\eta
_n,\mathcal{F}_n^j(\eta_n),L_2(P))}{\eta_n^2\sqrt{n}}M \biggr) \\[-2pt]
&\lesssim& q_n^{1/2}\eta_n + q_nn^{-1/2} \\[-2pt]
&=& O\bigl\{n^{\nu/2-\min(2\nu,(1-\nu)/2)}\bigr\} + O(n^{\nu-1/2}) \\[-2pt]
&=& O\bigl\{n^{-\min(3\nu/2,1/2-\nu)}\bigr\}+ O(n^{\nu-1/2}) = o(1),
\end{eqnarray*}
where the last equality holds because $0<\nu<1/2$.
Thus by Markov's inequality, $I_{1n} = n^{-1/2}\mathbb{G}_n \dot
{l}_\zeta(\hat{\theta}_n;Z)[h_j^*-h_{j,n}^*] = o_p(n^{-1/2})$.

Next for $I_{2n}$, the Taylor expansion for $\dot{l}_\zeta(\hat{\theta
}_n;Z)[h_j^*-h_{j,n}^*]$ at $\theta_0$ yields
\begin{eqnarray*}
&& \dot{l}_\zeta(\hat{\beta}_n,\hat{\zeta}_n(\cdot,\hat
{\beta}_n);Z)[h_j^*-h_{j,n}^*] - \dot{l}_\zeta(\beta_0,\zeta_0(\cdot
,\beta_0);Z)[h_j^*-h_{j,n}^*] \\[-2pt]
&&\qquad = (\hat{\beta}_n-\beta_0)'\ddot{l}_{\beta\zeta}(\tilde{\beta
}_n,\tilde{\zeta}_n(\cdot,\tilde{\beta}_n);Z)[h_j^*-h_{j,n}^*] \\[-2pt]
&&\qquad\quad{} + \ddot{l}_{\zeta\zeta}(\tilde{\beta}_n,\tilde{\zeta
}_n(\cdot,\tilde{\beta}_n);Z)[h_j^*-h_{j,n}^*,\hat{\zeta}_n-\zeta_0],
\end{eqnarray*}
where $(\tilde{\beta}_n,\tilde{\zeta}_n(\cdot,\tilde{\beta}_n))$ is
between $(\beta_0,\zeta_0(\cdot,\beta_0))$ and $(\hat{\beta}_n,\hat
{\zeta}_n(\cdot,\hat{\beta}_n))$. Then it follows that
\begin{eqnarray*}
&& |\ddot{l}_{\beta\zeta}(\tilde{\beta}_n,\tilde{\zeta
}_n(\cdot,\tilde{\beta}_n);Z)[h_j^*-h_{j,n}^*]| \\[-2pt]
&&\qquad = \biggl| X\biggl\{\Delta(\dot{w}_j^*-\dot{w}_{j,n}^*)(\epsilon
_{\tilde{\beta}_n})\\[-2pt]
&&\qquad\quad\hspace*{16pt}{} -\int_a^b 1(\epsilon_0 \geq t) \exp\{\tilde
{g}_n(t_{\tilde{\beta}_n})\}[(\dot{w}_j^*-\dot
{w}_{j,n}^*)(t_{\tilde{\beta}_n}) \\[-2pt]
&&\hspace*{144pt}\qquad\quad{} + \dot{\tilde{g}}_n(t_{\tilde{\beta}_n})
(w_j^*-w_{j,n}^*)(t_{\tilde{\beta}_n})] \,dt \biggr\} \biggr| \\[-2pt]
&&\qquad \lesssim\|\dot{w}_j^*-\dot{w}_{j,n}^* \|_{\infty} + \|\dot
{w}_j^*-\dot{w}_{j,n}^* \|_{\infty} \biggl\{\int_a^b \exp\{\tilde
{g}_n(t_{\tilde{\beta}_n})\} \,dt \biggr\} \\[-2pt]
&&\qquad\quad{} + \| w_j^*-w_{j,n}^* \|_{\infty} \biggl\{ \int_a^b \exp\{
\tilde{g}_n(t_{\tilde{\beta}_n})\}\dot{\tilde{g}}_n(t_{\tilde{\beta
}_n}) \,dt \biggr\} \\[-2pt]
&&\qquad \lesssim\|\dot{w}_j^*-\dot{w}_{j,n}^* \|_{\infty} + \|
w_j^*-w_{j,n}^* \|_{\infty} \\[-2pt]
&&\qquad = O(n^{-\nu})+O(n^{-2\nu})\\[-2pt]
&&\qquad=O(n^{-\nu}),
\end{eqnarray*}
where the second inequality holds because $\tilde{g}_n$ and its first
derivative $\dot{\tilde{g}}_n$ are bounded (or growing with $n$ slowly
enough so it can be effectively treated as bounded based on the same
argument of \cite{ShenWong1994} on page 591), and the last equality
holds due to the Corollary 6.21 of \cite{Schumaker} that $\|\dot
{w}_j^*-\dot{w}_{j,n}^* \|_{\infty}=O(n^{-(2-1)\nu})=O(n^{-\nu})$. Thus,
\begin{eqnarray*}
&&P|(\hat{\beta}_n-\beta_0)'\ddot{l}_{\beta\zeta}(\tilde{\beta
}_n,\tilde{\zeta}_n(\cdot,\tilde{\beta}_n);Z)[h_j^*-h_{j,n}^*]|\\
&&\qquad= |\hat{\beta}_n-\beta_0| \cdot O(n^{-\nu}) \\
&&\qquad= O_p\bigl\{n^{-\min(p\nu,(1-\nu)/2)}\bigr\}\cdot O(n^{-\nu}) \\
&&\qquad= O_p\bigl\{
n^{-\min((p+1)\nu,(1+3\nu)/2)}\bigr\}.
\end{eqnarray*}
Also,
\begin{eqnarray*}
&& |\ddot{l}_{\zeta\zeta}(\tilde{\beta}_n,\tilde{\zeta}_n(\cdot,\tilde
{\beta}_n);Z)[h_j^*-h_{j,n}^*,\hat{\zeta}_n-\zeta_0]| \\
&&\qquad = \biggl| \int_a^b 1(\epsilon_0 \geq t)\exp\{\tilde
{g}_n(t_{\tilde{\beta}_n})\} (w_j^*-w_{j,n}^*)(t_{\tilde{\beta}_n})(\hat
{g}_n-g_0)(t_{\tilde{\beta}_n}) \,dt \biggr| \\
&&\qquad \leq\| w_j^*-w_{j,n}^* \|_{\infty} \cdot\biggl\{ \int_a^b
\exp\{\tilde{g}_n(t_{\tilde{\beta}_n})\}(\hat{g}_n-g_0)(t_{\tilde{\beta
}_n}) \,dt\biggr\} \\
&&\qquad = \| w_j^*-w_{j,n}^* \|_{\infty}\cdot I_{3n}.
\end{eqnarray*}
By the Cauchy--Schwarz inequality and the boundedness of $\tilde{g}_n$,
we have
\begin{eqnarray*}
P\{I_{3n}\}^2 &=& P\biggl\{ \int_a^b \exp\{\tilde{g}_n(t_{\tilde{\beta
}_n})\}(\hat{g}_n-g_0)(t_{\tilde{\beta}_n}) \,dt\biggr\}^2 \\
&\lesssim& \int_{\mathcal{X}}\int_a^b (\hat{g}_n-g_0)^2(t_{\tilde{\beta
}}) \,d \Lambda_0(t) \,d F_X(x)
= \| \hat{\zeta}_n(\cdot,\tilde{\beta}_n) - \zeta_0(\cdot,\tilde{\beta
}_n) \|_2^2 \\
&\lesssim& |\tilde{\beta}_n-\hat{\beta}_n|^2 + \|\hat{\zeta}_n(\cdot
,\hat{\beta}_n)-\zeta_0(\cdot,\beta_0)\|_2^2 + |\beta_0-\tilde{\beta
}_n|^2 \\
&\lesssim& |\hat{\beta}_n-\beta_0|^2+ \| \hat{\zeta}_n(\cdot,\hat{\beta
}_n) - \zeta_0(\cdot,\beta_0) \|_2^2 = d(\hat{\theta}_n,\theta_0)^2.
\end{eqnarray*}
Hence $P|I_{3n}| \lesssim d(\hat{\theta}_n,\theta_0)$ and
\begin{eqnarray*}
&& P|\ddot{l}_{\zeta\zeta}(\tilde{\beta}_n,\tilde
{g}_n;Z)[h_j^*-h_{j,n}^*,\hat{\zeta}_n-\zeta_0]|\\
&&\qquad \lesssim\| w_j^*-w_{j,n}^*\|_{\infty} \cdot d(\hat{\theta
}_n,\theta_0)
= O(n^{-2\nu})\cdot O_p\bigl\{n^{-\min(p\nu,(1-\nu)/2)}\bigr\}\\
&&\qquad = O_p\bigl\{
n^{-\min((p+2)\nu,(1+3\nu)/2)}\bigr\}.
\end{eqnarray*}
Since $\frac{1}{2(1+p)} < \nu< \frac{1}{1+2p}$, it follows that
$ I_{2n} = O\{n^{-\min((p+1)\nu,(1+3\nu)/2)}\} = o(n^{-1/2})$.
Thus $I_n = I_{1n}+I_{2n} = o_p(n^{-1/2})$, and condition (A4) holds.

Now we\vspace*{1pt} verify assumption (A5). First by Lemma \ref{lemlem6},
the $\varepsilon$-bracketing numbers for the classes of functions
$\mathcal{F}_{n,j}^{\beta}(\eta)$
and
$\mathcal{F}_{n,j}^{\zeta}(\eta)$
are both bounded by $(\eta/\varepsilon)^{cq_n+d}$, which implies
that\vadjust{\goodbreak}
the corresponding $\varepsilon$-bracketing integrals are both bounded
by $q_n^{1/2}\eta$, that is,
\[
J_{[\csp]}(\eta,\mathcal{F}_{n,j}^{\beta}(\eta),L_2(P)) \lesssim
q_n^{1/2}\eta\quad\mbox{and}\quad J_{[\csp]}(\eta,\mathcal{F}_{n,j}^{\zeta
}(\eta),L_2(P)) \lesssim q_n^{1/2}\eta.
\]
Then for $\dot{l}_{\beta_j}(\theta;z)-\dot{l}_{\beta_j}(\theta_0;z)$,
by applying the Cauchy--Schwarz inequality, together with subtracting
and adding the terms $\dot{g}(\epsilon_0)$, $e^{g_0(t_{\beta})}\dot
{g}(t_{\beta})$, $e^{g_0(t)}\dot{g}(t_{\beta})$ and $e^{g_0(t)}\dot
{g}_0(t_{\beta})$, we have
\begin{eqnarray*}
&& \{\dot{l}_{\beta_j}(\theta;Z)-\dot{l}_{\beta
_j}(\theta_0;Z)\}^2 \\
&&\qquad= \biggl\{-\Delta
X_j[\dot{g}(\epsilon_\beta)-\dot{g}_0(\epsilon_0)]\\
&&\qquad\quad\hspace*{4pt}{}+X_j \int_a^b 1(\epsilon_0 \geq t)\bigl[e^{g(t_\beta)}\dot{g}(t_\beta
)-e^{g_0(t)}\dot{g}_0(t)\bigr] \,dt\biggr\}^2 \\
&&\qquad \lesssim\{\Delta[\dot{g}(\epsilon_\beta)-\dot{g}_0(\epsilon_0)]^2\}
+ \biggl\{
\int_a^b \bigl[e^{g(t_\beta)}\dot{g}(t_\beta)-e^{g_0(t)}\dot{g}_0(t)\bigr]^2
\,dt\biggr\} \\
&&\qquad \lesssim\{\Delta[\dot{g}(\epsilon_\beta)-\dot{g}(\epsilon_0)]^2\} +
\{\Delta[\dot{g}(\epsilon_0)-\dot{g}_0(\epsilon_0)]^2\} \\
&&\qquad\quad{} + \int_a^b \bigl\{ \bigl[e^{g(t_\beta)}-e^{g_0(t_\beta)}\bigr]^2 +
\bigl[e^{g_0(t_\beta)}-e^{g_0(t)}\bigr]^2 \bigr\} \dot{g}^2(t_\beta) \,dt \\
&&\qquad\quad{} + \int_a^b e^{2g_0(t)}\bigl\{[\dot{g}(t_\beta)-\dot
{g}_0(t_\beta)]^2 + e^{2g_0(t)}[\dot{g}_0(t_\beta)-\dot{g}_0(t)]^2
\bigr\} \,dt \\
&&\qquad = B_1+B_2+B_3+B_4.
\end{eqnarray*}

For $B_1$, since $\ddot{g}$ is bounded and the largest eigenvalue of
$P(XX')$ satisfies $0<\lambda_d<\infty$ by condition (C.2)(b), it follows that
\begin{eqnarray*}
PB_1 &\leq& P[\ddot{g}(Y-X'\tilde{\beta})X'(\beta-\beta_0)]^2 \lesssim
P[X'(\beta-\beta_0)]^2 \\
&\leq&\lambda_d |\beta-\beta_0|^2 \lesssim|\beta-\beta_0|^2
\leq\eta^2.
\end{eqnarray*}
For $B_2$, we have
\begin{eqnarray*}
PB_2 &\leq& \int_{\mathcal{X}} \biggl\{\int_a^b \bigl(\dot g(t)- \dot
g_0(t)\bigr)^2 \,d\Lambda_0(t) \biggr\}\,d F_X(x) \\
&=& \| \dot{g}(\psi(\cdot,\beta_0))-\dot{g}_{0}(\psi(\cdot,\beta_0)) \|
_2^2 \\
&\lesssim& |\beta- \beta_0|^2+ \| \dot{g}(\psi(\cdot,\beta))-\dot
{g}_{0}(\psi(\cdot,\beta_0)) \|_2^2 \lesssim\eta^2.
\end{eqnarray*}
For $B_3$, by using the mean value theorem, it follows that
\begin{eqnarray*}
PB_3 &=& P\biggl\{\int_a^b \bigl\{\bigl[e^{\tilde{g}(t_\beta)}(g-g_0)(t_\beta
)\bigr]^2 + \bigl[e^{g_0(t_{\tilde{\beta}})}X'(\beta-\beta_0)\bigr]^2 \bigr\} \dot
{g}^2(t_\beta) \,dt\biggr\} \\[-2pt]
&\lesssim& \int_{\mathcal{X}}\int_a^b (g-g_0)^2(t_\beta) \,d \Lambda
_0(t) \,d F_X(x) + P[X'(\beta-\beta_0)]^2 \\[-2pt]
&\lesssim& \| \zeta(\cdot,\beta)-\zeta_0(\cdot,\beta_0) \|_2^2 + |\beta
-\beta_0|^2 \leq\eta^2,
\end{eqnarray*}
where $\tilde{g} = g_0 + \xi(g-g_0)$ for some $0<\xi<1$ and thus is
bounded. Finally for~$B_4$, by the mean value theorem, it follows that
\begin{eqnarray*}
PB_4 &=& P \biggl\{\int_a^b e^{2g_0(t)}\bigl\{[\dot{g}(t_\beta)-\dot
{g}_0(t_\beta)]^2 + e^{2g_0(t)}[\dot{g}_0(t_\beta)-\dot{g}_0(t)]^2
\bigr\} \,dt \biggr\} \\[-2pt]
&\lesssim& \int_{\mathcal{X}}\int_a^b (\dot{g}-\dot{g}_0)^2(t_\beta) \,d
\Lambda_0(t) \,d F_X(x) + P \int_a^b [\ddot{g}_0(t_{\tilde{\beta
}})X'(\beta-\beta_0)]^2 \,dt \\[-2pt]
&\lesssim& \| \dot{g}(\psi(\cdot,\beta))-\dot{g}_0(\psi(\cdot,\beta)) \|
_2^2 + P[X'(\beta-\beta_0)]^2 \\[-2pt]
&\lesssim& \| \dot{g}(\psi(\cdot,\beta))-\dot{g}_0(\psi(\cdot,\beta_0))
\|_2^2 + |\beta-\beta_0|^2 \lesssim\eta^2.
\end{eqnarray*}
Therefore\vspace*{1pt} we have $ P\{\dot{l}_{\beta_j}(\theta;Z)-\dot{l}_{\beta
_j}(\theta_0;Z)\}^2
\lesssim\eta^2$.
Using the similar argument, we can show that $P\{\dot{l}_\zeta
(\theta;Z)[h_j^*]-\dot{l}_\zeta(\theta_0;Z)[h_j^*]\}^2 \lesssim
\eta^2$. By Lemma \ref{lemlem1}, we also have $\| \dot{l}_{\beta
_j}(\theta;Z)-\dot{l}_{\beta_j}(\theta_0;Z) \|_{\infty}$ and $\| \dot
{l}_{\zeta}(\theta;Z)[h_j^*]-\dot{l}_{\zeta}(\theta_0;Z)[h_j^*] \|
_{\infty}$ are both bounded. Now we pick $\eta$ as $\eta_n = O\{n^{-\min
((p-1)\nu,(1-\nu)/2)}\}$, then by the maximal inequality in Lemma 3.4.2
of \cite{vanWellner1996}, it follows that
\begin{eqnarray*}
E_P\| \mathbb{G}_n \|_{\mathcal{F}_{n,j}^{\beta}(\eta_n)}
&\lesssim& q_n^{1/2}\eta_n + q_nn^{-1/2} \\[-2pt]
&=& O\bigl\{n^{\max(({3}/{2}-p)\nu,\nu-{1}/{2})}\bigr\}+
O(n^{\nu-{1}/{2}})
= o(1),
\end{eqnarray*}
where the last equality holds since $p\geq3$ and $\nu< \frac{1}{2}$.
Similarly, we have
$E_P\| \mathbb{G}_n \|_{\mathcal{F}_{n,j}^\zeta(\eta_n)}= o(1)$. Thus
for $\xi=\min(p\nu,(1-\nu)/2)$ and
\[
Cn^{-\xi} = O\bigl\{n^{-\min(p\nu,(1-\nu
)/2)}\bigr\}
\]
by Markov's inequality,
\begin{eqnarray*}
\sup_{d(\theta,\theta_0) \leq Cn^{-\xi}}\mathbb{G}_n\{\dot{l}_{\beta
_j}(\beta,\zeta(\cdot,\beta);Z)-\dot{l}_{\beta_j}(\beta_0,\zeta_0(\cdot
,\beta_0);Z)\} &=& o_p(1),
\\[-2pt]
\sup_{d(\theta,\theta_0) \leq Cn^{-\xi}}\mathbb{G}_n\{\dot{l}_{\zeta
}(\beta,\zeta(\cdot,\beta);Z)[h_j^*]-\dot{l}_{\zeta}(\beta_0,\zeta
_0(\cdot,\beta_0);Z)[h_j^*]\} &=& o_p(1).
\end{eqnarray*}
This completes the verification of assumption (A5).

Finally, assumption (A6) can be verified by using the Taylor expansion.
Since the proofs for the two equations in (A6) are essentially
identical, we just prove the first equation. In a neighborhood of
$\theta_0\dvtx\{\theta\dvtx d(\theta,\theta_0) \leq Cn^{-\xi}, \theta\in
\Theta
_n^p\}$ with $\xi= \min(p\nu,(1-\nu)/2)$, the Taylor expansion for
$\dot{l}_{\beta}(\theta;Z)$ yields
\begin{eqnarray*}
\dot{l}_{\beta}(\theta;Z)
&=& \dot{l}_{\beta}(\theta_0;Z) + \ddot{l}_{\beta\beta}(\tilde{\theta
};Z)(\beta-\beta_0)+\ddot{l}_{\beta g}(\tilde{\theta};Z)[\zeta(\cdot
,\beta)-\zeta_0(\cdot,\beta_0)] \\[-2pt]
&=& \dot{l}_{\beta}(\theta_0;Z) + \ddot{l}_{\beta\beta}(\theta
_0;Z)(\beta-\beta_0) + \ddot{l}_{\beta\zeta}(\theta_0;Z)[\zeta(\cdot
,\beta)-\zeta_0(\cdot,\beta_0)] \\[-2pt]
&&{} + \{\ddot{l}_{\beta\beta}(\tilde{\theta};Z)(\beta-\beta
_0)-\ddot{l}_{\beta\beta}(\theta_0;Z)(\beta-\beta_0)\} \\[-2pt]
&&{} + \{\ddot{l}_{\beta\zeta}(\tilde{\theta};Z)[\zeta(\cdot
,\beta)-\zeta_0(\cdot,\beta_0)]-\ddot{l}_{\beta\zeta}(\theta_0;Z)[\zeta
(\cdot,\beta)-\zeta_0(\cdot,\beta_0)]\},
\end{eqnarray*}
where $\tilde{\theta}=(\tilde{\beta},\tilde{\zeta}(\cdot,\tilde{\beta
}))$ is a midpoint between $\theta_0$ and $\theta$. So
\begin{eqnarray*}
&& P\{\dot{l}_{\beta}(\theta;Z)-\dot{l}_{\beta}(\theta_0;Z) - \ddot
{l}_{\beta\beta}(\theta_0;Z)(\beta-\beta_0)-\ddot{l}_{\beta\zeta
}(\theta_0;Z)[\zeta(\cdot,\beta)-\zeta_0(\cdot,\beta_0)]\} \\
&&\qquad = P\{\ddot{l}_{\beta\beta}(\tilde{\theta};Z)-\ddot{l}_{\beta\beta
}(\theta_0;Z)\}(\beta-\beta_0) \\
&&\qquad\quad{} + P\{\ddot{l}_{\beta\zeta}(\tilde{\theta};Z)[\zeta(\cdot
,\beta)-\zeta_0(\cdot,\beta_0)]-\ddot{l}_{\beta\zeta}(\theta_0;Z)[\zeta
(\cdot,\beta)-\zeta_0(\cdot,\beta_0)]\}.
\end{eqnarray*}
Then by direct calculation we have
\begin{eqnarray*}
&& P|\ddot{l}_{\beta\beta}(\tilde{\theta};Z)-\ddot
{l}_{\beta\beta}(\theta_0;Z)| \\
&&\qquad\leq P|XX'\Delta\{\ddot{\tilde{g}}(\epsilon_{\tilde{\beta}})-\ddot
{g}_0(\epsilon_0)\}| \\
&&\qquad\quad{} + P\biggl\{XX'\biggl|\int_a^b 1(\epsilon_0 \geq t)\bigl\{\exp
\{\tilde{g}(t_{\tilde{\beta}})\} \ddot{\tilde{g}}(t_{\tilde{\beta
}})-\exp\{g_0(t)\}\ddot{g}_0(t)\bigr\} \,dt \\
&&\qquad\quad\hspace*{49.4pt}{} + \int_a^b 1(\epsilon_0 \geq t)\bigl\{\exp\{\tilde
{g}(t_{\tilde{\beta}})\} \dot{\tilde{g}}{}^2(t_{\tilde{\beta}})-\exp\{
g_0(t)\}\dot{g}^2_0(t)\bigr\} \,dt \biggr| \biggr\} \\
&&\qquad\lesssim P|\Delta\{\ddot{\tilde{g}}(\epsilon_{\tilde{\beta}})-\ddot
{g}_0(\epsilon_0)\}| \\
&&\qquad\quad{} + P\biggl\{\int_a^b |\exp\{\tilde{g}(t_{\tilde{\beta}})\}
\ddot{\tilde{g}}(t_{\tilde{\beta}})-\exp\{g_0(t)\}\ddot{g}_0(t)| \,dt
\biggr\} \\
&&\qquad\quad{} + P\biggl\{\int_a^b |\exp\{\tilde{g}(t_{\tilde{\beta}})\}
\dot{\tilde{g}}{}^2(t_{\tilde{\beta}})-\exp\{g_0(t)\}\dot{g}^2_0(t)| \,dt
\biggr\} \\
&&\qquad = C_1+C_2+C_3.
\end{eqnarray*}
By applying a similar argument that we used before for verifying (A5)
and condition (C.6), we can show
\begin{eqnarray*}
C_1 &\lesssim&|\beta-\beta_0| + \| \ddot{g}(\psi(\cdot,\beta))-\ddot
{g}_0(\psi(\cdot,\beta_0)) \|_2 \\
&=& O(n^{-\xi}) + O\bigl\{n^{-\min((p-2)\nu,
(1-\nu)/2)}\bigr\}.
\end{eqnarray*}
Similarly, we can show
\begin{eqnarray*}
C_2 &\lesssim&|\beta-\beta_0| + \| \ddot{g}(\psi(\cdot,\beta))-\ddot
{g}_0(\psi(\cdot,\beta_0)) \|_2\\
&=& O(n^{-\xi}) + O\bigl\{n^{-\min((p-2)\nu,
(1-\nu)/2)}\bigr\}
\end{eqnarray*}
and
\begin{eqnarray*}
C_3 &\lesssim&|\beta-\beta_0| + \| \dot{g}(\psi(\cdot,\beta))-\dot
{g}_0(\psi(\cdot,\beta_0)) \|_2 \\
&=& O(n^{-\xi}) + O\bigl\{n^{-\min((p-1)\nu,
(1-\nu)/2)}\bigr\},
\end{eqnarray*}
where $\xi=\min(p\nu, (1-\nu)/2)$. Therefore,
\[
P|\ddot{l}_{\beta\beta}(\tilde{\theta};Z)-\ddot{l}_{\beta\beta}(\theta
_0;Z)|= O\bigl\{n^{-\min((p-2)\nu, (1-\nu)/2)}\bigr\}
\]
and thus
\begin{eqnarray*}
&& P|\ddot{l}_{\beta\beta}(\tilde{\theta};Z)-\ddot{l}_{\beta\beta
}(\theta_0;Z)|(\beta-\beta_0) \\
&&\qquad = O\bigl\{n^{-\min((p-2)\nu, (1-\nu)/2)}\bigr\}\cdot O\bigl\{n^{-\min(p\nu,
(1-\nu)/2)}\bigr\} \\
&&\qquad = O\bigl\{n^{-\min(2(p-1)\nu,{1}/{2}+(p-{5}/{2})\nu,1-\nu
)}\bigr\} \\
&&\qquad= o(n^{-1/2}),
\end{eqnarray*}
where the last equality holds since $p \geq3$, so $2(p-1)\nu>\frac
{p-1}{p+1} \geq\frac{1}{2}$, $\frac{1}{2}+(p-\frac{5}{2})\nu>\frac
{1}{2}$ and $1-\nu>\frac{1}{2}$. Similarly we can show
\begin{eqnarray*}
&& P|\ddot{l}_{\beta\zeta}(\tilde{\theta};Z)[\zeta(\cdot,\beta)-\zeta
_0(\cdot,\beta_0)]-\ddot{l}_{\beta\zeta}(\theta_0;Z)[\zeta(\cdot,\beta
)-\zeta_0(\cdot,\beta_0)]| \\
&&\qquad = O\bigl\{n^{-\min(2(p-1)\nu,{1}/{2}+(p-{5}/{2})\nu,1-\nu
)}\bigr\} \\
&&\qquad= o(n^{-1/2}).
\end{eqnarray*}
Therefore, we have
\begin{eqnarray*}
&& |P\{\dot{l}_{\beta}(\theta;Z)-\dot{l}_{\beta}(\theta_0;Z) - \ddot
{l}_{\beta\beta}(\theta_0;Z)(\beta-\beta_0)-\ddot{l}_{\beta\zeta
}(\tilde{\theta};Z)[\zeta(\cdot,\beta)-\zeta_0(\cdot,\beta_0)]\}| \\
&&\qquad = O\bigl\{n^{-\min(2(p-1)\nu,{1}/{2}+(p-{5}/{2})\nu,1-\nu)}\bigr\}
=
O(n^{-\alpha\xi}),
\end{eqnarray*}
where
$\alpha= \min(2(p-1)\nu,\frac{1}{2}+(p-\frac{5}{2})\nu,1-\nu)/\min(p\nu
,\frac{1-\nu}{2}) > 1 $ and $\alpha\xi> 1/2$.\vspace*{1pt}

Therefore, we have verified all six assumptions, and thus we have
\[
\sqrt{n}(\hat{\beta}_n-\beta_0)=A^{-1}\sqrt{n}\mathbb{P}_nl_{\beta
_0}^*(\beta_0,\zeta_0(\cdot,\beta_0);Z) + o_p(1) \rightarrow
N(0,A^{-1}B(A^{-1})'),
\]
where $l_{\beta_0}^*(\theta_0;Z)=\dot{l}_{\beta}(\theta_0;Z)-\dot
{l}_\zeta(\theta_0;Z)[\mathbf{h}^*] $ is the efficient score function
for~$\beta_0$ and $A=P\{{l}^*_{\beta_0}(Y, \Delta, X)\}^{\otimes
2}=I(\beta_0)$, which is shown when verifying (A3). Hence $A=B$ and
$A^{-1}B(A^{-1})'=A^{-1}=I^{-1}(\beta_0)$, and
\[
\sqrt{n}\mathbb{P}_n l_{\beta_0}^*(\theta_0;Z) = n^{-{1}/{2}}\sum
_{i=1}^n l^*_{\beta_0}(Y_i,\Delta_i,X_i).
\]
Thus we complete the proof of Theorem \ref{thmnormality}.

\section*{Acknowledgments}

The authors would like to thank two referees and an associate editor
for their very helpful comments.

\begin{supplement}
\stitle{Additional proofs}
\slink[doi]{10.1214/11-AOS934SUPP} 
\sdatatype{.pdf}
\sfilename{aos934\_supp.pdf}
\sdescription{The supplementary document contains proofs of technical
lemmas and Theorems \ref{thmConvRate} and \ref{thmvariance}.}
\end{supplement}

%

\printaddresses

\end{document}